\documentclass[a4paper,11pt]{article}
\usepackage{ifstyle}
\usepackage{graphicx}
\usepackage{epsfig} 

\usepackage[francais]{babel}               
\usepackage[T1]{fontenc}
\usepackage[latin1]{inputenc}
\usepackage{amssymb,diagrams}             
\usepackage{array}                 
\usepackage{fullpage}
\newtheorem{theoreme}{Théorème} 
\newtheorem{lem}[subsubsection]{Lemme}

\newcommand{\N}{\mathbb N}

\newcommand{\R}{\mathbb R}
\newcommand{\C}{\mathbb C}
\newcommand{\K}{\mathbb K}
\newcommand{\I}{\mathbb I}

\begin{document}
      
\begin{center} 
\huge{\bf{Une caractérisation différentielle des faisceaux analytiques cohérents sur une variété complexe}}
\Large{Nefton Pali}
\end{center} 
{\bf Résumé}.-Dans ce travail nous généralisons, dans le contexte des faisceaux, un résultat classique de Grothendieck concernant l'intégrabilité des connexions de type $(0,1)$ sur un fibré vectoriel ${\cal C}^{\infty}$ au dessus d'une variété complexe. En introduisant la notion de faisceau $\bar{\partial}$-cohérent, qui est une notion qui vit dans le contexte ${\cal C}^{\infty}$, nous montrons l'existence d'une équivalence (exacte) entre la catégorie des faisceaux analytiques cohérents et la catégorie des faisceaux $\bar{\partial}$-cohérents. La difficulté essentielle de la preuve de ce résultat consiste à résoudre un problème différentiel quasi-linéaire dont le terme principal est un opérateur $\bar{\partial}$ usuel. La solution de ce problème est obtenue en utilisant un schéma de convergence rapide de type Nash-Moser. Un autre ingrédient important pour la conclusion de la preuve est un résultat profond de Malgrange qui affirme la fidélité plate de l'anneau des germes des fonctions ${\cal C}^{\infty}$ à valeurs complexes en un point, sur l'anneau des germes des fonctions holomorphes en ce point.
\\
\\
{\bf{Abstract}}.-In this work we give a generalization, in the context of sheaves, of a classical result of Grothendieck concerning the integrability of connections of type $(0,1)$ over a ${\cal C}^{\infty}$ vector bundle over a complex manifold. We introduce the notion of $\bar{\partial}$-coherent sheaf, which is a ${\cal C}^{\infty}$ notion, and we prove the existence of an (exact) equivalence between the category of coherent analytic sheaves and the category of $\bar{\partial}$-coherent sheaves. The principal difficulty of the proof is the solution of a quasi-linear differential equation with standard $\bar{\partial}$ as its principal term. We are able to find a solution of this differential equation, using a rapidly convergent iteration scheme of Nash-Moser type. We also use a deep result of Malgrange asserting that the ring of germs of complex differentiable functions at a point is faithfully flat over the ring of germs of holomorphic functions at the same point.
\section{Introduction }
Le résultat de Grothendieck mentionné dans le résumé (voir \cite{ko-mal}) affirme qu'un fibré vectoriel complexe différentiable au dessus d'une variété complexe, qui admet une connexion $\bar{\partial}$ de type $(0,1)$ telle que $\bar{\partial}^2=0$, possède une structure de fibré vectoriel holomorphe. En autres termes, en utilisant l'équivalence entre les notions de fibré vectoriel et de faisceau localement libre sur une variété,  on a que le noyau de la connexion  $\bar{\partial}$ est un faisceau de ${\cal O}$-modules localement libre. Le point essentiel de ce résultat consiste a prouver l'existence d'une solution de l'équation différentielle quasi-linéaire  $g^{-1}\bar{\partial}_{_J }g=A$ (où $\bar{\partial}_{_J }$ est la (0,1)-connexion canonique sur le faisceau des fonctions ${\cal C}^{\infty}$ et $A$ représente la $(0,1)$-forme de connexion locale de $\bar{\partial}$) avec la condition d'intégrabilité  $\bar{\partial}_{_J }A +A\wedge A=0$. Nous généraliserons cette équation pour prouver notre caractérisation différentielle laquelle introduit comme objet nouveau la notion de faisceau $\bar{\partial}$-cohérent. Précisément un faisceau  $\bar{\partial}$-cohérent est un faisceau ${\cal G}$ de modules de fonctions ${\cal C}^{\infty}$ à valeurs complexes, muni d'une connexion $\bar{\partial}$ de type $(0,1)$ telle que $\bar{\partial}^2=0$, et qui admet localement  des résolutions de longueur finie, par des  modules de fonctions ${\cal C}^{\infty}$ à valeurs complexes. Notre caractérisation affirme essentiellement que le noyau de la connexion $\bar{\partial}$ est un faisceau analytique cohérent. On aura alors une équivalence exacte (l'exactitude est due à la fidélité plate de l'anneau des germes des fonctions ${\cal C}^{\infty}$ à valeurs complexes sur l'anneau des germes des fonctions holomorphes, voir \cite{mal}) entre la catégorie des faisceaux analytiques cohérents et la catégorie des faisceaux  $\bar{\partial}$-cohérents. La difficulté essentielle de la preuve consiste à montrer que quel que soit le choix de la résolution locale du faisceau ${\cal G}$, on peut trouver au voisinage de chaque point de l'ouvert sur lequel on considère la résolution locale, une autre résolution locale constituée de matrices holomorphes. En d'autres termes on cherche une résolution locale  admettant des formes de connexion nulles. Pour atteindre cet'objectif on introduit la notion de recalibration, laquelle généralise la notion classique de changement de jauge pour les formes locales d'une connexion de type $(0,1)$ intégrable sur un faisceau localement libre. La notion de recalibration ne représente rien d'autre que une action d'un semi-groupe sur l'ensemble des formes qui représentent localement la condition d'intégrabilité de la connexion $\bar{\partial}$. La notion en question permet de traduire notre problème en terme d'un système différentiel quasi-linéaire dont terme principal un opérateur $\bar{\partial}$ usuel. Les conditions d'intégrabilité de ce système ne sont rien d'autre que les expressions locales de la condition d'intégrabilité de la connexion $\bar{\partial}$. Les solutions du système différentiel seront obtenues à l'aide d'un procédé itératif de type Nash-Moser, dont chaque étape  est déterminée par une recalibration obtenue en fonction de l'étape précédente. La preuve de l'existence de solutions du système différentiel en question est exposée dans la sub-section (3.4), qui constitue la partie essentielle de la preuve de notre caractérisation des faisceaux analytiques cohérents. La technique qui consiste à utiliser des schémas itératifs pour montrer l'existence des solutions de problèmes non linéaires est désormais bien connue en analyse complexe. On cite par exemple les travaux de Webster \cite{We-1} et \cite{We-2}, qui utilise des techniques de type Nash-Moser, (voir \cite{Mos-1} et \cite{Mos-2}) pour prouver l'existence des solutions de deux problèmes différentiels fondamentaux en géométrie complexe. L'ingrédient final qui permet de conclure notre preuve est le  résultat profond de Malgrange cité précédemment, lequel permet aussi une généralisation du théorème de Dolbeault au cas des faisceaux analytiques cohérents ($\bar{\partial}$-cohérents).  Enfin on remarque aussi un résultat d'intégrabilité pour les connexions des faisceaux admettant des résolutions locales de longueur finie sur les variétés différentiables.

%%%%%%%%%%%%%%%%%%%%%%%%%%%%%%%%%%%%%%%%%%%%%%%%%%%%%%%%%%%%%%%%%%%%%%%%%
\section{Faisceaux $\bar{\partial}$-cohérents sur les variétés complexes}
%%%%%%%%%%%%%%%%%%%%%%%%%%%%%%%%%%%%%%%%%%%%%%%%%%%%%%%%%%%%%%%%%%%%%%%%%

Soit $(X,J)$ une variété complexe, où $J$ est le tenseur, supposé  intégrable, de la structure presque-complexe. On désigne par ${\cal E}_X$ le faisceau
 des fonctions ${\cal C}^{\infty} $ à valeurs complexes, par  ${\cal E}^{0,q}_X$ le faisceau des $(0,q)$-formes et par
$\bar{\partial}_{_J }\in Hom_{_{{\cal O}_X } }({\cal E}^{0,q}_X, {\cal E}^{0,q+1}_X )(X)$ la composante de type $(0,1)$ de la différentielle. Sur tout faisceau $\cal F$ de ${\cal O}_X$-modules on peut considérer
 la connexion canonique 
$$\bar{\partial}_{_{\cal F} }:=\I_{_{\cal F}} \otimes _{_{{\cal O}_X }}\bar{\partial}_{_J }
:{\cal F}\otimes_{_{{\cal O}_X }}{\cal E}^{0,q}_X 
\longrightarrow  {\cal F} \otimes_{_{{\cal O}_X }}{\cal E}^{0,q+1}_X $$
qui est vue comme une connexion de type $(0,1)$ sur le faisceau de  ${\cal E}_X$-modules ${\cal F}^{\infty}:={\cal F} \otimes_{_{{\cal O}_X }}{\cal E}_X$. De manière générale, sur un faisceau $\cal G$ de
 ${\cal E}_X$-modules il suffit de donner un morphisme de faisceaux de groupes additifs
 $\bar{\partial}:{\cal G}\longrightarrow  {\cal G} \otimes_{_{{\cal E}_X }}{\cal E}^{0,1}_X $
 tel 
que $\bar{\partial}(g\cdot f)=( \bar{\partial} g)\cdot f +g\otimes\bar{\partial}_{_J } f$, $g\in {\cal G}_x$, $f\in {\cal E}_{X,x}$, pour déterminer de façon univoque une connexion $\bar{\partial}$ de type $(0,1)$ sur le complexe $({\cal G} \otimes_{_{{\cal E}_X }}{\cal E}^{0,q}_X)_{ q \geq 0}$. En effet on peut définir l'extension $\bar{\partial}:{\cal G}\otimes_{_{{\cal E}_X }}{\cal E}^{0,q}_X\longrightarrow  {\cal G} \otimes_{_{{\cal E}_X }}{\cal E}^{0,q+1}_X $ par la formule classique
\begin{eqnarray*} 
(\bar{\partial}\omega )(\xi _0,...,\xi _q):=\sum_{0\leq j \leq q}(-1)^j \bar{\partial}(\omega (\xi _0,...,\widehat{\xi _j},..., \xi _q))(\xi _j)+
\\
+\sum_{0\leq j<k \leq q}(-1)^{j+k}\omega ([\xi _j,\xi _k],\xi _0,...,\widehat{\xi _j},...,\widehat{\xi _k},..., \xi _q)  
\end{eqnarray*} 
avec $\omega\in ({\cal G} \otimes_{_{{\cal E}_X }}{\cal E}^{0,q}_X)(U)$ et $\xi _j\in{\cal E}(T^{0,1}_X)(U)$, sur un ouvert $U$ quelconque, ou de façon équivalente, par la règle de Leibniz
 $\bar{\partial}(g\otimes\alpha ):=\bar{\partial}g \wedge \alpha +g\otimes\bar{\partial}_{_J }\alpha  $ avec $g\in{\cal G}_x$ et $\alpha \in {\cal E}^{0,q}_{X,x} $ pour tout $x\in X$. La donnée d'une connexion de type $(0,1)$ sur ${\cal G}$ détermine aussi le tenseur de 
courbure de la connexion qu'on 
notera $\Theta _{\bar{\partial}}\in ({\cal E}nd_{_{{\cal E}_{_X}}}({\cal G}) \otimes_{_{{\cal E}_{_X} }}{\cal E}^{0,2}_X)(X)$ et que on défini par la formule $\Theta _{\bar{\partial}}(\xi ,\eta)\cdot g:=(\bar{\partial}^2g)(\xi ,\eta)$ 
avec $g\in{\cal G}(U)$ et $\xi ,\eta\in{\cal E}(T^{0,1}_X)(U)$. On note de plus $\xi _{\bar{\partial} }\cdot g:=\bar{\partial}g(\xi )$ la dérivée covariante de la section $g$ calculée le long du champ $\xi $ et on remarque que la première définition de l'extension de $\bar{\partial}$ implique de façon triviale la formule
$$\xi _{\bar{\partial} }\cdot(\eta _{\bar{\partial} }\cdot g)- \eta _{\bar{\partial} }\cdot(\xi _{\bar{\partial} }\cdot g )=[\xi ,\eta]_{\bar{\partial} }\cdot g +\Theta _{\bar{\partial}}(\xi ,\eta)\cdot g
$$
où $[\xi ,\eta]\in{\cal E}(T^{0,1}_X)(U)$, grâce à l'hypothèse d'intégrabilité du tenseur de la structure presque-complexe $J\in {\cal C}^{\infty} (T_X^*\otimes_{_{_\R}} T_X)(X)$. Le tenseur de courbure $\Theta _{\bar{\partial}}$ exprime donc le défaut de commutation des dérivées covariantes secondes des sections de ${\cal G}$ le long des champs de type $(0,1)$. Nous porterons un intérêt particulier aux connexions de type $(0,1)$ intégrables, c'est à dire aux connexions telles que $\bar{\partial}^2=0$. La formule précédente caractérise alors ce type de connexions $(0,1)$ comme étant celles pour lesquelles les dérivées covariantes secondes, calculées le long de deux champs qui commutent, commutent egalément. En termes explicites on a l'égalité $\xi _{\bar{\partial} }\cdot(\eta _{\bar{\partial} }\cdot g)= \eta _{\bar{\partial} }\cdot(\xi _{\bar{\partial} }\cdot g )$ si $[\xi ,\eta]=0$. Un exemple de $(0,1)$-connexion intégrable est évidement la connexion $\bar{\partial}_{_{\cal F}}$ introduite précédemment.
Avant d'énoncer le résultat qu'on se propose de démontrer, on remarque que si $\cal F$ est un faisceau analytique cohérent, on dispose d'un diagramme commutatif dont toutes les directions horizontales et verticales sont exactes.    
\begin{diagram}[height=1cm,width=1cm]
&         &0&                 &0&                                         &0&                                             &0&                                    \\
&         &\uTo&              &\uTo&                                      &\uTo&                                          &\uTo&                                 \\
0&  \rTo  & {\cal F}
          _{|_U}&  \rTo   &{\cal F}^{\infty}_{|_U}&    \rTo^{\bar{\partial} _{_{\cal F}}} 
 & {\cal F}^{\infty}_{|_U}\otimes_{_{{\cal E}_{_U}  }}{\cal E}^{0,1}_{_U} &   \rTo^{\bar{\partial}_{_{\cal F}}} & {\cal F}^{\infty}_{|_U}\otimes_{_{{\cal E}_{_U}   }}{\cal E}^{0,2}_{_U} &\rTo& \cdot \cdot \cdot
\\
&        &\uTo_{\psi} &              &\uTo_{\psi} &                                      &\uTo_{\psi \otimes \I_{(0,1)} } &                                          &\uTo_{\psi \otimes \I_{(0,2)} } &                                 \\
0&\rTo  & {\cal O}
          _{_U}^{\oplus p_0} &  \rTo   &{\cal E}^{\oplus p_0}_{_U}&    \rTo^{\bar{\partial}_{_J }}
 & ({\cal E}^{0,1}_{_U} )^{\oplus p_0}&   \rTo^{\bar{\partial}_{_J }} & ({\cal E}^{0,2}_{_U} )^{\oplus p_0}&\rTo&\cdot \cdot \cdot
\\
&        &\uTo_{\varphi _1}  &              &\uTo_{\varphi _1} &                                      &\uTo_{\varphi _1\otimes \I_{(0,1)} } &                                          &\uTo_{\varphi _1\otimes \I_{(0,2)} } &                                 \\
0&\rTo  & {\cal O}
          _{_U}^{\oplus p_1} &  \rTo   &{\cal E}^{\oplus p_1}_{_U}&    \rTo^{\bar{\partial}_{_J }}
 & ({\cal E}^{0,1}_{_U} )^{\oplus p_1}&   \rTo^{\bar{\partial}_{_J }} & ({\cal E}^{0,2}_{_U} )^{\oplus p_1}&\rTo&\cdot \cdot \cdot
\\
&        &\uTo_{\varphi _2}  &              &\uTo_{\varphi _2} &                                      &\uTo_{\varphi _2\otimes \I_{(0,1)} } &                                          &\uTo_{\varphi _2\otimes \I_{(0,2)} } &                                                  
\\
&         &\vdots&                 & \vdots&                                         &\vdots&                                             &\vdots&             
\\
&        &\uTo_{\varphi_{m-1} }  &              &\uTo_{\varphi_{m-1} } &                                      &\uTo_{\varphi_{m-1} \otimes \I_{(0,1)} } &                                          &\uTo_{\varphi_{m-1} \otimes \I_{(0,2)} } &                                                  
\\
0&\rTo  & {\cal O}
          _{_U}^{\oplus p_{m-1} } &  \rTo   &{\cal E}^{\oplus p_ {m-1}}_{_U}&    \rTo^{\bar{\partial}_{_J }}
 & ({\cal E}^{0,1}_{_U} )^{\oplus p_{m-1} }&   \rTo^{\bar{\partial}_{_J }} & ({\cal E}^{0,2}_{_U} )^{\oplus p_{m-1} }&\rTo&\cdot \cdot \cdot
\\
&        &\uTo_{\varphi _m}  &              &\uTo_{\varphi _m} &                                      &\uTo_{\varphi _m\otimes \I_{(0,1)} } &                                          &\uTo_{\varphi _m\otimes \I_{(0,2)} } &                                 \\
0&\rTo  & {\cal O}
          _{_U}^{\oplus p_m} &  \rTo   &{\cal E}^{\oplus p_m}_{_U}&    \rTo^{\bar{\partial}_{_J }}
 & ({\cal E}^{0,1}_{_U} )^{\oplus p_m}&   \rTo^{\bar{\partial}_{_J }} & ({\cal E}^{0,2}_{_U} )^{\oplus p_m}&\rTo&\cdot \cdot \cdot
\\
&        &\uTo&              &\uTo&                                      &\uTo&                                          &\uTo&                                 \\
&         &0&                 &0&                                         &0&                                             &0&                                    \\       
\end{diagram}
\\ 
La raison de l'exactitude est la suivante. Le théorème des syzygies implique qu'on peut choisir une ${\cal O}$-résolution de longueur finie du faisceau analytique cohérent ${\cal F}$ comme celle du diagramme précédent. Ensuite la platitude de l'anneau ${\cal E}_{_{X,x} }$ sur l'anneau ${\cal O}_{_{X,x} }$ (voir l'ouvrage de Malgrange \cite{mal}) implique l'exactitude des autres flèches verticales. L'exactitude du dernier complexe $(({\cal E}^{0,q}_{_U} )^{\oplus p_m};\bar{\partial}_{_J } )_{ q \geq 0}$ implique l'exactitude du complexe  
$({\cal R}^{\cal E}(\varphi _{m-1})   \otimes_{_{{\cal E}_U }}{\cal E}^{0,q}_U;\bar{\partial}_{_J } )_{ q \geq 0}$ où ${\cal R}^{\cal E}(\varphi _{m-1})$ désigne le faisceau des relations de $\varphi _{m-1}$ sur le faisceau ${\cal E}_{_{X} }  $. En procédent par récurrence décroissante et en utilisant l'exactitude des complexes en $\bar{\partial}_{_J }$ et l'exactitude des flèches verticales on obtient finalement l'exactitude du complexe 
$({\cal F}^{\infty}_{|_U}\otimes_{_{{\cal E}_{_U}  }}{\cal E}^{0,q}_{_U};\bar{\partial} _{_{\cal F}})_{ q \geq 0} $. De manière générale on a la caractérisation différentielle suivante.
\begin{theoreme}
Soit X une variété complexe et soit ${\cal G}$ un faisceau de  ${\cal E}_X$-modules,  muni d'une connexion 
$\bar{\partial}:{\cal G}\longrightarrow  {\cal G} \otimes_{_{{\cal E}_X }}{\cal E}^{0,1}_X $ de type $(0,1)$ telle que $\bar{\partial}^{2}=0$. Si de plus le faisceau ${\cal G}$ admet localement une ${\cal E}$-résolution de longueur finie, alors le faisceau de  ${\cal O}_X$-modules $Ker\bar{\partial}\subset\cal G$ est analytique cohérent, on a les égalités 
${\cal G}=(Ker\bar{\partial})\cdot {\cal E}_X\cong
(Ker\bar{\partial})\otimes_{_{{\cal O}_{_X}}}{\cal E}_X$ et  la connexion $\bar{\partial}$ coïncide, à isomorphisme canonique prés, avec l'extension naturelle $\bar{\partial}_{Ker\bar{\partial}} $ associée au faisceau analytique cohérent $Ker\bar{\partial} $.
\end{theoreme}
 Considérons maintenant la définition suivante.
\begin{defi}
Un couple $({\cal G},\bar{\partial})\equiv{\cal G}_{\bar{\partial}}$ où ${\cal G}$ et $\bar{\partial}$ vérifient les hypothèses du théorème 1 est appelé faisceau $\bar{\partial}$-cohérent. Un morphisme $\varphi :{\cal A}_{\bar{\partial}_1}\longrightarrow {\cal B}_{\bar{\partial}_2}$ de faisceaux $\bar{\partial}$-cohérents est un 
morphisme de faisceaux de ${\cal E}_X$-modules telles que le diagramme suivant soit commutatif
\begin{diagram}[height=1cm,width=1cm]
{\cal A}\otimes_{_{{\cal E}_{_X}}}{\cal E}^{0,1}_X&\rTo^{\varphi\otimes \I_{(0,1)} }&{\cal B}\otimes_{_{{\cal E}_{_X}}}{\cal E}^{0,1}_X
\\
\uTo^{\bar{\partial}_1}&              &\uTo^{\bar{\partial}_2}  
\\
{\cal A}&\rTo^{\varphi}&\cal B
\end{diagram}  
\end{defi}
Dans le cas des faisceaux de ${\cal E}_X$-modules inversibles  qui admettent une connexion $\bar{\partial}_0$ de type $(0,1)$ telle que $\bar{\partial}^{2}_0=0$ on sait que toutes les connexions de ce type, et seulement celles ci, sont de la forme  $\bar{\partial}_0+A\otimes$ où $A\in{\cal E}^{0,1}_X(X)$ est une $(0,1)$ -forme $\bar{\partial}_{_J }$-fermée.
\\
Le théorème 1 et la fidélité plate du faisceau ${\cal E}_X$ sur ${\cal O}_X$ montrent que sur une variété complexe on a une équivalence exacte entre la catégorie ${\cal O}{\bf Coh} $ 
des faisceaux analytiques cohérents et la catégorie $\bar{\partial}{\bf Coh} $ des faisceaux $\bar{\partial}$-cohérents. Plus explicitement on a le foncteur $\infty$ qui agit de la façon suivante
\begin{eqnarray*}
{\cal F}\in{\cal O}{\bf Coh} &\stackrel{\infty}{\longmapsto}&{\cal F}^{\infty}_{\bar{\partial}_{_{\cal F}}}\in\bar{\partial}{\bf Coh} 
\\
 \varphi \in Hom_{_{{\cal O}_X}}({\cal A},{\cal B})&\stackrel{\infty}{ \longmapsto}&\varphi\otimes \I\in 
Hom({\cal A}^{\infty}_{\bar{\partial}_{_{\cal A}}},{\cal B}^{\infty}_{\bar{\partial}_{_{\cal B}}})
\\
\\
{\cal G}_{\bar{\partial}} \in\bar{\partial}{\bf Coh}  & \stackrel{\infty^{-1}}{ \longmapsto}&Ker\bar{\partial} \in{\cal O}{\bf Coh} 
\\
 \varphi \in Hom({\cal A}_{\bar{\partial}_1},{\cal B}_{\bar{\partial}_2}) &\stackrel{\infty^{-1}}{ \longmapsto}   
 & \varphi_{|..}\in Hom_{_{{\cal O}_X}}(Ker\bar{\partial}_1,Ker\bar{\partial}_2)
\end{eqnarray*}
Une conséquence immédiate de l'équivalence exacte précédente est l'existence d'une équivalence exacte entre la catégorie des faisceaux d'ideaux de fonctions holomorphes cohérents et la  catégorie des faisceaux d'ideaux ${\cal J}\subseteq {\cal E}_X$ de fonctions ${\cal C}^{\infty}$ à valeurs complexes admettant des ${\cal E}$-résolutions locales de longueur finie, qui sont stables par rapport aux dérivations le long des champs de vecteurs de type (0,1), autrement dit les faisceaux d'ideaux ${\cal J}\subseteq {\cal E}_X$ tels que pour tout germes de (0,1)-champs $\xi \in{\cal E}(T^{0,1}_X)_x$ on a l'inclusion $\xi .{\cal J}_x\subseteq {\cal J}_x $, pour tout point $x\in X$. En termes plus concrets, si on se donne un faisceau d'ideaux ${\cal J}\subseteq {\cal E}_X$ de fonctions ${\cal C}^{\infty}$ à valeurs complexes admettant des ${\cal E}$-résolutions locales de longueur finie, on a que le faisceau d'ideaux 
${\cal J}\cap {\cal O}_X$ est analytique cohérent si pour un choix arbitraire de repaires locales 
$(\xi_1,...,\xi_n)\in{\cal E}(T^{0,1}_X)^{\oplus n}(U)$ et de générateurs $(\psi_1,...,\psi_p)\in{\cal J}^{\oplus p}(U)$ on a pour tout $x\in U$, l'existence de germes de fonctions $f_{k,l,j}\in {\cal E}_{_{X,x} }$  qui vérifient les égalités 
$
\xi _{k,x} \,.\psi_{l,x}=\sum_{j=1}^p\,f_{k,l,j}\cdot\psi_{j,x}
$
. On remarque enfin que la cohomologie des faisceaux cohérents ($\bar{\partial}$-cohérents) sur une variété complexe peut se calculer, grâce à l'isomorphisme fonctoriel de De Rham-Weil (voir par exemple l'ouvrage de Demailly \cite{Dem})  par la formule suivante:
$$
H^q(X,{\cal G}_{\bar{\partial}}):=H^q (X,Ker\bar{\partial})\cong H^q(\Gamma(X,{\cal G} \otimes_{_{{\cal E}_X }}{\cal E}^{0,*}_X);\bar{\partial})
$$ 
qui constitue une généralisation du théorème de Dolbeault. Un cas particulier (ou une généralisation si on veut) de la formule précédente est la suivante: 
$$
H^q(X,{\cal G}_{\bar{\partial}}  \otimes_{_{{\cal E}_X }}{\cal E}^{p,0}_{\bar{\partial}_{_{J,p} }  }) :=H^q(X,(Ker\bar{\partial})\otimes_{_{{\cal O}_X }}{\cal O}(\Omega ^p_X)) \cong H^q(\Gamma(X,{\cal G} \otimes_{_{{\cal E}_X }}{\cal E}^{p,*}_X);\bar{\partial}_{\pi} )=:H^{p,q}(X,{\cal G}_{\bar{\partial}})$$
avec $\bar{\partial}_{_{J,p} }:=(-1)^p  \bar{\partial}_{_J }$, ${\cal G}_{\bar{\partial}}  \otimes_{_{{\cal E}_X }}{\cal E}^{p,0}_{\bar{\partial}_{_{J,p} } }:=({\cal G} \otimes_{_{{\cal E}_X }}
{\cal E}^{p,0}_X \,,\bar{\partial}_{\pi })$ et $\bar{\partial}_{\pi }:{\cal G}\otimes_{_{{\cal E}_X }}{\cal E}^{p,0}_X\longrightarrow  {\cal G} \otimes_{_{{\cal E}_X }}{\cal E}^{p,1}_X $ définie par la règle de Leibniz
.

%%%%%%%%%%%%%%%%%%%%%%%%%%%%%%%%%%%%%%%%%%%%%%%%%%%%%%%%%%%%%%%%%%%%%%%%%%%%%%%%%%%%%%%%%%%%%%%%%%%%
\section{Preuve du théorème 1}
\subsection{Première étape: expression locale de la condition d'intégrabilité $\bar{\partial}^2=0$.} 
%%%%%%%%%%%%%%%%%%%%%%%%%%%%%%%%%%%%%%%%%%%%%%%%%%%%%%%%%%%%%%%%%%%%%%%%%%%%%%%%%%%%%%%%%%%%%%%%%%%%

Nous commençons par prouver le lemme élémentaire suivant.
\begin{lem} 
Soit X une variété complexe et soit ${\cal G}$ un faisceau de  ${\cal E}_X$-modules.
\\
$(A)$ Supposons que le faisceau ${\cal G}$ admet des ${\cal E}$-présentations localement finies et soit
$$
 {\cal E}_{_U}^{\oplus p_1 }\stackrel{\varphi}{ \longrightarrow }{\cal E}_{_U}^{\oplus p_0 }\stackrel{\psi}{\longrightarrow}
{\cal G}_{_{|_U} }\rightarrow 0
$$
une ${\cal E}$-présentation au dessus d'un ouvert $U$. Alors l'existence d'une connexion 
$\bar{\partial}$ de type $(0,1)$ sur le faisceau ${\cal G}_{_{|_U} }$ telle que 
$\bar{\partial}^{2}=0$, implique l'existence de matrices $\omega ^{s,0}\in M_{p_s,p_s}({\cal E}^{0,1}_X(U)),\;s=0,1$ et 
$\omega ^{0,1}\in M_{p_1,p_0}({\cal E}^{0,2}_X(U))$ telles que $\bar{\partial}\psi=\psi\cdot \omega^{0,0}$ et
\begin{eqnarray}
\bar{\partial}_{_J }\varphi+\omega ^{0,0}\cdot \varphi =  \varphi\cdot \omega ^{1,0} \quad
\\\nonumber
\\
\bar{\partial}_{_J }\omega ^{0,0} +\omega ^{0,0}\wedge \omega ^{0,0}=\varphi \cdot\omega ^{0,1} 
\end{eqnarray}
Réciproquement l'existence de matrices $\omega ^{s,0},\;s=0,1 $ et $ \omega ^{0,1} $ qui vérifient les relations $(1)$ et $(2)$, implique l'existence d'une connexion $\bar{\partial}$ de type $(0,1)$ sur le faisceau ${\cal G}_{_{|_U} }$ telle que 
$\bar{\partial}\psi=\psi\cdot \omega^{0,0}$ et
$\bar{\partial}^{2}=0$. 
\\
$(B)$ Supposons que le faisceau ${\cal G}$ admet des ${\cal E}$-résolutions locales de longueur finie et soit 
$$
0\rightarrow {\cal E}_{_U}^{\oplus p_m }\stackrel{\varphi_m }{\longrightarrow}  
{\cal E}_{_U}^{\oplus p_{m-1} }\stackrel{ \varphi _{m-1} }{\longrightarrow} \cdot \cdot\cdot
\stackrel{\varphi_2}{ \longrightarrow } {\cal E}_{_U}^{\oplus p_1 }\stackrel{\varphi_1}{ \longrightarrow }{\cal E}_{_U}^{\oplus p_0 }\stackrel{\psi}{\longrightarrow}
{\cal G}_{_{|_U} }\rightarrow 0
$$
une telle ${\cal E}$-résolution. Alors l'existence d'une connexion 
$\bar{\partial}$ de type $(0,1)$ sur le faisceau ${\cal G}_{_{|_U} }$ telle que 
$\bar{\partial}^{2}=0$, implique l'existence des matrices $\omega^{s,k}  \in M_{p_{s+k} ,p_s}({\cal E}^{0,k+1}_X(U))$ pour $s=0,...,m$ et $k=-1,...,m-s$ telles que si on utilise l'identification $\varphi _s\equiv\omega ^{s,-1}$ et les conventions formelles $\omega ^{0,-1} :=0,\;\omega ^{-1,j}:=0$ et $\omega ^{s,k}:=0$ si $s\geq m+1$ ou $k\geq m-s+1$, on aura les relations $\bar{\partial}\psi=\psi\cdot \omega^{0,0}$ et
\begin{eqnarray}
\bar{\partial}_{_J }\omega^{s,k} +\sum_{j=-1}^{k+1} (-1)^{k-j}  \omega^{s+j,k-j}\wedge \omega^{s,j}=0 
\end{eqnarray}
pour tout $s=0,...,m$, et $k=-1,...,m-s,\;(s,k)\not= (0,-1)$.
\\
Réciproquement l'existence des matrices $\omega ^{s,k}$ qui vérifient la relation $(3)$ implique l'existence d'une connexion  $\bar{\partial}$ de type $(0,1)$ sur le faisceau ${\cal G}_{_{|_U} }$ telle que 
$\bar{\partial}\psi=\psi\cdot \omega^{0,0}$ et
$\bar{\partial}^{2}=0$. 
\end{lem} 
$Preuve\; de\;(A) $. La ${\cal E}$-présentation de ${\cal G}_{_{|_U}}$ considérée dans l'hypothèse implique l'existence des ${\cal E}$-présentations
$$
({\cal E}^{0,q}_{_U})^{\oplus p_1}\longrightarrow ({\cal E}^{0,q}_{_U})^{\oplus p_0}\longrightarrow  {\cal G}_{_{|_U}}\otimes_{_{{\cal E}_{_U}}}{\cal E}^{0,1}_{_U}\longrightarrow 0
$$ 
pour $q\geq 0$. On aura alors l'existence d'une matrice $\omega ^{0,0}\in M_{p_0,p_0}({\cal E}^{0,1}_X(U))$ telle que $\bar{\partial}\psi=\psi\cdot \omega^{0,0}$. En appliquant la connexion $\bar{\partial}$ à l'identité $\psi \circ \varphi=0$ on obtient $\psi(\bar{\partial}_{_J }\varphi +\omega ^{0,0}\cdot  \varphi)=0$. L'exactitude des 
${\cal E}$-présentations précédentes, implique alors l'existence d'une matrice 
$\omega ^{1,0}\in M_{p_1,p_1}({\cal E}^{0,1}_X(U))$ telle que la relation (1) soit satisfaite. On obtient alors le diagramme commutatif suivant, ayant des flèches verticales exactes:

\begin{diagram}[height=1cm,width=1cm]
0&                                       &0&                                                                &0
\\
\uTo&                                    &\uTo&                                                            &\uTo
\\
{\cal G}_{_{|_U} }&\rTo^{\bar{\partial}} &{\cal G}_{_{|_U}}\otimes_{_{{\cal E}_{_U}}}{\cal E}^{0,1}_{_U}
                                                                               &\rTo^{\bar{\partial}} &    {\cal G}_{_{|_U}}
                                                                                                        \otimes_{_{{\cal E}_{_U}}}{\cal E}^{0,2}_{_U} 
\\
\uTo^{\psi}&                             &\uTo_{\psi \otimes \I_{(0,1)} }&                               &\uTo_{\psi \otimes \I_{(0,2)} }
\\
{\cal E}^{\oplus p_0} _{_U}&\rTo^{\bar{\partial}_{_J }+\omega ^{0,0} }&({\cal E}^{0,1}_{_U})^{\oplus p_0} &\rTo^{\bar{\partial}_{_J }+\omega ^{0,0} }&({\cal E}^{0,2}_{_U})^{\oplus p_0} 
\\
\uTo^{\varphi_1}&                             &\uTo_{\varphi_1\otimes \I_{(0,1)} }&                               &\uTo_{\varphi_1 \otimes \I_{(0,2)} }
\\
{\cal E}^{\oplus p_1} _{_U}&\rTo^{\bar{\partial}_{_J }+\omega ^{1,0} }&({\cal E}^{0,1}_{_U})^{\oplus p_1}&\rTo^{\bar{\partial}_{_J }+\omega ^{1,0} }&({\cal E}^{0,2}_{_U})^{\oplus p_1} 
\end{diagram} 
\\
\\
L'hypothèse d'intégrabilité $\bar{\partial}^2=0$ implique $\psi(\bar{\partial}_{_J }\omega ^{0,0} +\omega ^{0,0}  \wedge \omega ^{0,0} )=0$, d'où l'existence d'une matrice $\omega^{0,1}  \in M_{p_1,p_0}({\cal E}^{0,2}_X(U))$ telle que la relation (2) soit satisfaite.
\\
Pour prouver la réciproque dans la partie (A), il suffit de considérer la connexion quotient $\bar{\partial} $ obtenue par la connexion $\bar{\partial}_{_J }+\omega ^{0,0}\wedge \bullet $.
\\
$Preuve\;de\;(B)$. Pour $(s,k)=(1,-1)$ et $(s,k)=(0,0)$ la relation (3) exprime les relations (1) et (2) prouvées dans la partie (A) de la preuve du lemme. En  applicant l'opérateur $\bar{\partial}_{_J }$ aux identités $\varphi_{t-1} \circ \varphi_t=0$ on obtient inductivement, de la même façon utilisée pour obtenir la relation (1), l'existence d'une matrice $\omega ^{s,0} \in M_{p_s,p_s}({\cal E}^{0,1}_X(U))$, $s=1,...,m$ telle que 
$
\bar{\partial}_{_J }\varphi _s +\omega ^{t-1,0} \cdot  \varphi_s =  \varphi_s\cdot \omega ^{t,0}  
$
. Ces relations constituent les relations (3) pour $(s,-1),\;s=1,...,m$. On va montrer maintenant l'existence des matrices $\omega ^{s,k},\;k\geq 1 $ qui vérifient la relation (3) à l'aide du procédé récursif triangulaire suivant. Pour un couple $(s,k)$, $s=0,...,m-1$ et $k=1,...,m-s$ on suppose avoir déjà défini $\omega^{\sigma ,\kappa}$ pour $\sigma +\kappa \leq s+k$, $\kappa \leq k+1$, et on applique l'opérateur $\bar{\partial}_{_J }$ à l'expression $(3_{s,k-1})$ pour $k\geq 1$. On obtient alors la relation suivante:
\begin{eqnarray*}
 \sum_{j=-1}^{k} (-1)^{k-j-1}\bar{\partial}_{_J }  \omega^{s+j,k-j-1}\wedge \omega^{s,j}- \sum_{j=-1}^{k}  \omega^{s+j,k-j-1}\wedge \bar{\partial}_{_J } \omega^{s,j} =0 
\end{eqnarray*}
En explicitant les termes $\bar{\partial}_{_J }  \omega^{\bullet,\bullet}$ dans la relation précédente (qui bien évidemment, grâce à l'hypothèse de récurence précédente, vérifient $(3_{\bullet,\bullet})$ pour les indices voulus) on obtient:
\begin{eqnarray*}
\sum_{j=-1}^{k} (-1)^{k-j-1}  \omega^{s+k,-1}\wedge \omega^{s+j,k-j}\wedge \omega^{s,j} + \qquad\qquad\qquad
\\
\\
+\sum_{j=-1}^{k} \sum_{r=-1}^{k-j-1}(-1)^{r+1} 
\omega^{s+j+r,k-j-1-r}\wedge \omega^{s+j,r} \wedge \omega^{s,j}+\qquad\qquad
\\
\\
+ \sum_{j=-1}^{k-1} \sum_{r=-1}^{j+1}(-1)^{j-r} 
\omega^{s+j,k-j-1}\wedge \omega^{s+r,j-r} \wedge \omega^{s,r}-\omega^{s+k,-1}\wedge\bar{\partial}_{_J }  \omega^{s,k}  =0 
\end{eqnarray*}
En faisant le changement d'indice $j'=j+r$, $r'=j$ dans la deuxième somme et en rappelant que $\omega^{s-1,-1} \wedge \omega^{s,-1}=0$ on obtient: 
\begin{eqnarray*}
\omega^{s+k,-1}\wedge\big(\bar{\partial}_{_J }\omega^{s,k}+\sum_{j=-1}^{k} (-1)^{k-j} \omega^{s+j,k-j}\wedge \omega^{s,j}\big)=0  
\end{eqnarray*}
L'hypothèse d'exactitude  nous permet de choisir $\omega^{s,k+1}$ telle que la relation
\begin{eqnarray*}
\bar{\partial}_{_J }\omega^{s,k}+\sum_{j=-1}^{k} (-1)^{k-j} \omega^{s+j,k-j}\wedge \omega^{s,j}=\omega^{s+k+1,-1}\wedge\omega^{s,k+1}  
\end{eqnarray*}
soit satisfaite.
 Ce type de récurrence peut se visualiser grâce au tableau suivant:
%%%%% ICI
\begin{figure}[hbtp]
\begin{center} 
\input dgt0.pstex_t       
   \caption{}
    \label{fig1}
\end{center} 
\end{figure}    
%%%%%%%%
\\
L'hypothèse de finitude de la longueur des ${\cal E}$-résolutions locales de ${\cal G}$ permet d'arrêter  ce procédé  après un nombre fini d'étapes. On a donc prouvé la première implication de la partie (B) du lemme. Le réciproque dans la partie (B) est évidemment une conséquence banale de la partie (A) du lemme.
\\
La figure 2 montre les  matrices $\omega^{s,k}$  qui sont représentées par des flèches dans le diagramme suivant lequel représente le complexe déterminé par la 
${\cal E}$-résolution locale $(\varphi ,\psi)$ dans le cas de longueur $m=4$.
\\
\\
%%%%% ICI
\begin{figure}[hbtp]
\begin{center} 
\input dgt1.pstex_t
    \caption{}
    \label{fig2}
\end{center} 
\end{figure}    
%%%%%%%% 
\\
Le diagramme suivant montre les matrices qui interviennent dans la relation $(3)$ pour $(s,k)=(1,2)$, dans le cas de longueur $m=4$.
%%%%% ICI
\begin{figure}[hbtp]
\begin{center} 
\input dgt2.pstex_t       
    \caption{}
    \label{fig3}
\end{center} 
\end{figure}    
%%%%%%%%
\\
La liberté homologique qui caractérise le choix des matrices $\omega ^{s,k}$ est exprimée par une action de semi-groupe qui aura une importance considérable dans la preuve du théorème 1 et que on expose dans la deuxième partie de la preuve.
\\
{\bf{Remarque.} } A partir de maintenant le lecteur doit tenir compte du fait que certains des calculs et formules qui suivront n'existent pas dans le cas de longueur $m=0$ de la ${\cal E}$-résolutions locale, autrement dit dans le cas des faisceaux localement libres. Cependant les calculs qui survivent ont encore  sens et ils font partie de notre preuve (différente de la preuve donnée par Grothendieck) dans ce cas. On utilisera aussi la convention qui consiste à négliger les termes d'une somme ou d'un produit si l'ensemble des indices sur lesquels on effectue ces opérations est vide. 

%%%%%%%%%%%%%%%%%%%%%%%%%%%%%%%%%%%%%%%%%%%%%%%%%%%%%%%%%%%%%%%
\subsection{Deuxième étape: le formalisme du procédé itératif.}
%%%%%%%%%%%%%%%%%%%%%%%%%%%%%%%%%%%%%%%%%%%%%%%%%%%%%%%%%%%%%%%
 
Dans cette partie on se propose de présenter un formalisme utile pour la suite. On commence avec les définitions suivantes. On pose par définition
$$
\Gamma(U):={\displaystyle  \bigoplus_{\scriptstyle  s=0,...,m} }GL(p_s,{\cal E}_X(U))
$$
\begin{defi}
La classe $[\varphi,\psi]$ de ${\cal E}$-isomorphisme de la ${\cal E}$-résolution locale $(\varphi ,\psi)$ est l'ensemble des ${\cal E}$-résolutions locales $(\tilde{\varphi} ,\tilde{\psi})$ de longueur $m$ au dessus de l'ouvert $U$ pour lesquelles il existe $g\in\Gamma(U)$ tel que le diagramme suivant soit commutatif.
\begin{diagram}[height=1cm,width=1cm]
0&\rTo &{\cal E}^{\oplus p_m}_{_U} &\rTo^{\tilde{\varphi}_m } &\cdot\cdot\cdot&\rTo^{\tilde{\varphi}_2}  &{\cal E}^{\oplus p_1}_{_U}&\rTo^{\tilde{\varphi}_1}   &{\cal E}_{_U}^{\oplus p_0}&\rTo^{\tilde{\psi}} & {\cal G}_{_{|_U} }&\rTo &0
\\
&       &\dTo^{g_m}_{\wr}&                                         &          &                                &\dTo^{g_1}_{\wr}&          &\dTo^{g_0}_{\wr}&   &\dTo_{\I} 
\\
0&\rTo &{\cal E}^{\oplus p_m}_{_U}  &\rTo^{\varphi _m}       &\cdot\cdot\cdot&\rTo^{\varphi _2}       &{\cal E}^{\oplus p_1}_{_U}  &\rTo^{\varphi _1}  &{\cal E}_{_U}^{\oplus p_0}&\rTo^{\psi}&{\cal G}_{_{|_U} }&\rTo & 0 
\end{diagram} 
\end{defi}
On a alors que $[\varphi,\psi]=\{(\varphi _g,\psi_g)\,|\,g\in\Gamma (U)\} $ où $\psi_g:=\psi\cdot g_0$ et $\varphi _{s,g}:=g_{s-1}^{-1}\cdot \varphi _s \cdot g_s$. Ensuite on désigne par
\begin{eqnarray*} 
\Omega (U,\varphi_g,\psi_g,\bar{\partial})\subset 
{\displaystyle  \bigoplus_
{{\scriptstyle  s=0,...,m    \atop   \scriptstyle  k=-1,...,m-s}
\atop
\scriptstyle (s,k)\not=(0,-1)} }
M_{p_{s+k} ,p_s}({\cal E}^{0,k+1}_X(U))
\end{eqnarray*} 
l'ensemble, non vide par l'hypothèse d'exactitude de la ${\cal E} $-résolution locale $(\varphi ,\psi)$, constitué par les éléments 
$\omega =(\omega ^{s,k})_{s,k},\;\omega ^{s,k}\in M_{p_{s+k} ,p_s}({\cal E}^{0,k+1}_X(U))$  tels que 
$ \omega ^{\bullet,-1}=\varphi _{\bullet,g},\;\bar{\partial}\psi_g =\psi_g\cdot \omega ^{0,0}$ et la relation (3) soit satisfaite. Ensuite on définit la ``fibration''
$$
\Omega (U,[\varphi,\psi],\bar{\partial}):= {\displaystyle  \coprod_{g\in\Gamma(U) }  }\Omega (U,\varphi_g,\psi_g,\bar{\partial})
$$
au dessus du groupe $\Gamma(U) $ et l'ensemble des paramètres au dessous de l'ouvert $U$
$$
{\cal P}(U):={\displaystyle  \bigoplus_{\scriptstyle  s=0,...,m\atop\scriptstyle  k=0,...,m-s} }M_{p_{s+k} ,p_s}({\cal E}^{0,k}_X(U))
$$
constitué par les éléments $\eta=(\eta^{s,k})_{s,k}$, $ \eta^{s,k}\in M_{p_{s+k} ,p_s}({\cal E}^{0,k}_X(U))$, tels que 
$\I_{p_s}+\eta^{s,0}\in GL(p_s,{\cal E}_X(U))$, $s=0,...,m$. On munit ${\cal P}(U)$ de la loi de semi-groupe donnée par le produit extérieur des paramètres qu'on définit de la façon suivante; si $\eta_1,\;\eta_2\in{\cal P}(U) $ on désigne par $\eta_1\wedge\eta_2\in{\cal P}(U) $ le paramètre dont les composantes sont définies par la formule   
$$
(\eta_1\wedge \eta_2)^{s,k}:= \eta_1^{s,k}+\eta_2^{s,k}+ \displaystyle{\sum_{j=0}^k\eta^{s+j,k-j}_1\wedge \eta^{s,j}_2 } 
$$
L'élément neutre de cette loi est le zéro.  Considérons maintenant l'application
$$
\begin{array}{ccc}
R :{\cal P}(U)\times\Omega (U,[\varphi,\psi],\bar{\partial})& \longrightarrow &\Omega (U,[\varphi,\psi],\bar{\partial})
\\
\\
(\eta,\omega )&\longmapsto &\omega _{\eta}  
\end{array}
$$   
 où les composantes de $\omega _{\eta}$ sont définies par récurrence  sur $k$, pour tout les indices $s=0,...,m$, $k=-1,...,m-s,\;(s,k)\not=(0,-1)$ par la formule   
\begin{eqnarray}
\omega^{s,k}_{\eta}= \Big(\I_{p_{s+k}} +\eta^{s+k,0}\Big)^{-1}\cdot \Big( \bar{\partial}_{_J }\eta^{s,k}+\sum_{j=0}^{k+1}  \omega^{s+j,k-j}\wedge \eta^{s,j}-\sum_{j=-1}^{k-1} (-1)^{k-j}\eta^{s+j,k-j}\wedge \omega^{s,j}_{\eta} +\omega^{s,k}\Big)  
\end{eqnarray}
Dans l'expression précédente on utilise les définitions formelles $\eta^{s,m-s+1}:=0$, $\eta^{-1,j}:=0$ et $\eta^{s,-1}:=0$. On appelle $R$ application de recalibration et on dit que $ \omega _{\eta} $ est la recalibration de $\omega $ avec paramètre de recalibration $\eta$. Dans la suite on utilisera aussi les notations $g_s\equiv g_s(\eta):=\I_{p_s}+\eta^{s,0}$, $\psi_{\eta}\equiv\psi_{g(\eta)}$. Avec la première notation on a évidement 
$\omega ^{\bullet,-1}_{\eta} =\omega ^{\bullet,-1}_{g(\eta)}$. Le diagramme suivant montre les matrices qui interviennent dans la définition de 
$\omega^{1,2}_{\eta}$ dans le cas où la longueur  de la résolution est égale à 4. 
%%%%% ICI
\begin{figure}[hbtp]
\begin{center} 
\input dgt3.pstex_t
    \caption{}
    \label{fig4}
\end{center} 
\end{figure}    
%%%%%%%%
\\
Le lecteur peut alors essayer d'avoir une perception visuelle de la formule de recalibration (4).         
On remarque aussi que si $\omega \in\Omega (U,\varphi ,\psi,\bar{\partial})$ alors pour tout $\eta \in {\cal P}(U)$, la recalibration 
$R(\eta,\omega )\equiv\omega _{\eta}$ de $\omega $ détermine le diagramme commutatif suivant:
\begin{diagram}[height=1cm,width=1cm]
0&                                       &0&                                                                &0
\\
\uTo&                                    &\uTo&                                                            &\uTo
\\
{\cal G}_{_{|_U} }&\rTo^{\bar{\partial}} &{\cal G}_{_{|_U}}\otimes_{_{{\cal E}_{_V}}}{\cal E}^{0,1}_{_U}
                                                                               &\rTo^{\bar{\partial}} &    {\cal G}_{_{|_U}}
                                                                                                        \otimes_{_{{\cal E}_{_U}}}{\cal E}^{0,2}_{_U} 
\\
\uTo^{\psi_{\eta} }&                             &\uTo_{\psi_{\eta}  \otimes \I_{(0,1)} }&                               &\uTo_{\psi_{\eta}  \otimes \I_{(0,2)} }
\\
{\cal E}^{\oplus p_0} _{_U}&\rTo^{\bar{\partial}_{_J }+\omega^{0,0}_{\eta} }&({\cal E}^{0,1}_{_U})^{\oplus p_0} &\rTo^{\bar{\partial}_{_J }+\omega^{0,0}_{\eta} }&({\cal E}^{0,2}_{_U})^{\oplus p_0} 
\\
\uTo^{\omega^{1,-1}_{\eta} }&                             &\uTo_{\omega^{1,-1}_{\eta} \otimes \I_{(0,1)} }&                               &\uTo_{\omega^{1,-1}_{\eta}  \otimes \I_{(0,2)} }
\\
{\cal E}^{\oplus p_1} _{_U}&\rTo^{\bar{\partial}_{_J }+\omega^{1,0}_{\eta}}&({\cal E}^{0,1}_{_U})^{\oplus p_1}&\rTo^{\bar{\partial}_{_J }+\omega^{1,0}_{\eta} }&({\cal E}^{0,2}_{_U})^{\oplus p_1} 
\\
\uTo&         &\uTo&          &\uTo
\\
\vdots&                                       &\vdots&                                                                &\vdots
\\
\end{diagram}
Si on désigne par ${\cal P}_0(U)\subset {\cal P}(U)$ le sous-ensemble des paramètres tels que $\eta^{\bullet,0}=0$ on a que la restriction de la recalibration 
$R_0 :{\cal P}_0(U)\times\Omega (U,[\varphi,\psi],\bar{\partial})\longrightarrow \Omega (U,[\varphi,\psi],\bar{\partial})$ est une application fibrée qui mesure la liberté homologique qui caractérise le choix des matrices $\omega^{\bullet,\bullet}$ relativement aux ${\cal E} $-résolutions locales $(\varphi _g,\psi_g)$, $g\in \Gamma(U)$. Avec les notations introduites précédemment on a la proposition suivante.
\begin{prop}
L'application de recalibration $R$ est bien définie et constitue une action de semi-groupe  transitive sur l'ensemble $\Omega (U,[\varphi,\psi],\bar{\partial})$
\end{prop}  
$Preuve$. Nous commençons par prouver que l'application $R $ est bien définie. 
On prouve d'abord les relations $\bar{\partial}\psi_{\eta} =\psi_{\eta}\cdot\omega^{0,0}_{\eta}$. En effet on a :
\begin{eqnarray*}
\bar{\partial}\psi_{\eta} =\bar{\partial}\psi\cdot g_0+ \psi \cdot \bar{\partial}\eta^{0,0} =
\psi\cdot(\bar{\partial}\eta^{0,0}+\omega ^{0,0}\wedge \eta^{0,0}+ \omega ^{0,0})=
\\
\\
=\psi\cdot(\bar{\partial}\eta^{0,0}+\omega ^{0,0}\wedge \eta^{0,0}+\omega ^{1,-1}\wedge \eta^{0,1} + \omega ^{0,0})=\psi_{\eta}\cdot \omega ^{0,0}_{\eta}\;\;
\end{eqnarray*}
On prouve maintenant que les matrices $\omega^{\bullet,\bullet}_{\eta} $ vérifient la relation (3) pour l'indice $k=-1$, autrement dit on veut montrer la relation:
$$
\bar{\partial}_{_J } \omega ^{s,-1}_{\eta} + \omega ^{s-1,0}_{\eta}\cdot \omega ^{s,-1}_{\eta} =\omega ^{s,-1}_{\eta}\cdot\omega ^{s,0}_{\eta}
$$
On commence par développer le terme $\bar{\partial}\omega ^{s,-1}_{\eta} $, en utilisant la relation (3) pour l'indice $k=-1$, relativement aux matrices $\omega^{\bullet,\bullet}$. On obtient alors les égalités suivantes:
\begin{eqnarray*}
\bar{\partial}_{_J } \omega ^{s,-1}_{\eta}=-g_{s-1}^{-1}\cdot\bar{\partial}_{_J } g_{s-1}\cdot g_{s-1}^{-1}\cdot\omega ^{s,-1} \cdot g_s+ \qquad\;\;
\\
\\
+g_{s-1}^{-1}\cdot  \bar{\partial}_{_J } \omega ^{s,-1}\cdot g_s+  g_{s-1}^{-1}\cdot   \omega ^{s,-1}\cdot \bar{\partial}_{_J }g_s=\qquad\quad
\\
\\
=-g_{s-1}^{-1}(\bar{\partial}_{_J }\eta^{s-1,0}+\omega ^{s-1,0}\wedge \eta^{s-1,0}+\omega ^{s-1,0})\cdot\omega ^{s,-1}_{\eta} +
\\
\\
+\omega ^{s,-1}_{\eta} \cdot g_s^{-1}(\bar{\partial}_{_J }\eta^{s,0}+\omega ^{s,0}\wedge \eta^{s,0}+\omega ^{s,0})\qquad \qquad
\end{eqnarray*}
En rappelant que $\omega ^{s-1,-1}\cdot \omega ^{s,-1}=0$ et en rajoutant et en soustrayant le terme 
$-g_{s-1}^{-1}\cdot \omega ^{s,-1}\cdot\eta^{s-1,1}\cdot \omega ^{s,-1}_{\eta}$ à la dernière expression de $\bar{\partial}\omega ^{s,-1}_{\eta} $, on obtient:
\begin{eqnarray*}
\bar{\partial}_{_J } \omega ^{s,-1}_{\eta}=-\omega ^{s-1,0}_{\eta}\cdot \omega ^{s,-1}_{\eta}+ \qquad \qquad \qquad 
\\
\\
+\omega ^{s,-1}_{\eta}\cdot g_s^{-1}(\bar{\partial}_{_J }\eta^{s,0}+\omega ^{s,0}\wedge \eta^{s,0}+\eta^{s-1,1}\wedge \omega ^{s,-1}_{\eta}+ \omega ^{s,0})= 
\\
\\
=-\omega ^{s-1,0}_{\eta}\cdot \omega ^{s,-1}_{\eta} + \omega ^{s,-1}_{\eta}\cdot\omega ^{s,0}_{\eta} \qquad \qquad\qquad
\end{eqnarray*}
On va montrer maintenant la validité de la formule (3) pour tous les indices, relativement aux matrices $\omega^{\bullet,\bullet}_{\eta}$, avec un procédé récursif analogue à celui qui nous a permis de définir les matrices $\omega^{\bullet,\bullet}$. Voici les détails de la récurrence. 
Pour un couple $(s,k)$, $s=0,...,m$, $k=0,...,m-s$ on suppose avoir déjà montré la relation
$$ 
\qquad \qquad \qquad \qquad \qquad \qquad \bar{\partial}_{_J }\omega^{\sigma ,\kappa}_{\eta} +\sum_{j=-1}^{\kappa+1} (-1)^{\kappa-j}  \omega^{\sigma +j,\kappa-j}_{\eta} \wedge \omega^{\sigma ,j}_{\eta}=0  \qquad \qquad \qquad \qquad \qquad  (3^{\sigma ,\kappa}_{\eta}) 
$$
pour $\sigma=s$, $\kappa=-1,...,k-1$. En développant le terme en $\bar{\partial}_{_J }$ de l'expression suivante on a l'égalité:
\begin{eqnarray*}
\bar{\partial}_{_J }\omega^{s,k}_{\eta} +\sum_{j=-1}^{k} (-1)^{k-j}  \omega^{s+j,k-j}_{\eta}\wedge \omega^{s,j}_{\eta}= \qquad\qquad\qquad\qquad\qquad
\\
\\
=g_{s+k}^{-1}\Big(-\bar{\partial}_{_J } \eta^{s+k,0}\wedge\omega^{s,k}_{\eta}+
\sum_{j=0}^{k+1}\bar{\partial}_{_J }\omega^{s+j,k-j}\wedge \eta^{s,j}- \sum_{j=0}^{k+1}(-1)^{k-j}\omega^{s+j,k-j}\wedge\bar{\partial}_{_J } \eta^{s,j}-
\\
\\
- \sum_{j=-1}^{k-1}(-1)^{k-j}\bar{\partial}_{_J }\eta^{s+j,k-j}\wedge\omega^{s,j}_{\eta}- \sum_{j=-1}^{k-1}\eta^{s+j,k-j}\wedge \bar{\partial}_{_J }\omega^{s,j}_{\eta}+\bar{\partial}_{_J }\omega^{s,k} \Big)+\qquad
\\
\\
+\sum_{j=-1}^{k} (-1)^{k-j}  \omega^{s+j,k-j}_{\eta}\wedge \omega^{s,j}_{\eta}=(A_1) \qquad\qquad\qquad\qquad\qquad
\end{eqnarray*}
En développant les termes $\bar{\partial}_{_J }\omega^{s+j,k-j}$ et $\bar{\partial}_{_J }\omega^{s,j}_{\eta}$ à l'aide respectivement des expressions $(3)$ et $(3^{s ,j}_{\eta})$ on obtient:
\begin{eqnarray*}
(A_1)=g_{s+k}^{-1}  \omega^{s+k+1,-1}\Big(\bar{\partial}_{_J } \eta^{s,k+1} +\sum_{j=0}^{k+1}\omega^{s+j,k+1-j}\wedge\eta^{s,j}\Big)+\qquad\qquad\qquad\quad
\\
\\
+g_{s+k}^{-1} \Big(\sum_{j=0}^{k+1}\sum_{r=-1}^{k-j}(-1)^{k-j-r+1}\omega^{s+j+r,k-j-r}\wedge\omega^{s+j,r}\wedge\eta^{s,j}- \sum_{j=0}^{k}(-1)^{k-j}\omega^{s+j,k-j}\wedge\bar{\partial}_{_J } \eta^{s,j}-
\\
\\
-\sum_{j=-1}^{k}(-1)^{k-j}\bar{\partial}_{_J }\eta^{s+j,k-j}\wedge\omega^{s,j}_{\eta}    + \sum_{j=-1}^{k-1}\sum_{r=-1}^{j+1}(-1)^{j-r}\eta^{s+j,k-j}\wedge \omega^{s+r,j-r}_{\eta}\wedge \omega^{s,r}_{\eta} +\bar{\partial}_{_J }\omega^{s,k}\Big)+
\\
\\
 +\sum_{j=-1}^{k} (-1)^{k-j}  \omega^{s+j,k-j}_{\eta}\wedge \omega^{s,j}_{\eta}=(A_2) \qquad\qquad\qquad\qquad\qquad\quad
\end{eqnarray*}
En faisant le changement d'indice $j'=j+r$, $r'=j$ dans la première somme double et en développant les premiers facteurs $\omega^{s+j,k-j}_{\eta}$ de la dernière somme on a:
\begin{eqnarray*}
(A_2)=g_{s+k}^{-1}\omega^{s+k+1,-1}\Big(\bar{\partial}_{_J } \eta^{s,k+1} +\sum_{j=0}^{k+1}\omega^{s+j,k+1-j}\wedge\eta^{s,j}\Big)+
g_{s+k}^{-1}\bar{\partial}_{_J }\omega^{s,k} 
\\
\\
+\sum_{j=-1}^{k}(-1)^{k-j+1}g_{s+k}^{-1}\omega^{s+j,k-j}\wedge\Big(\bar{\partial}_{_J } \eta^{s,j} + \sum_{r=0}^{j+1}\omega^{s+r,j-r}\wedge\eta^{s,r}\Big)+
\\
\\
+\sum_{j=-1}^{k} (-1)^{k-j}g_{s+k}^{-1} \Big(\sum_{r=0}^{k-j+1}\omega^{s+j+r,k-j-r}\wedge\eta^{s+j,r}+\omega^{s+j,k-j}\Big)\wedge \omega^{s,j}_{\eta}=(A_3)
\end{eqnarray*}
En rappelant que $g_{s+j,0}:=\I_{s+k}+\eta^{s+j,0}$, en décomposant les termes extrêmes de la somme $\sum_{r=0}^{k-j+1}$ et en décomposant les facteurs $\omega^{s,j}_{\eta}$ qui apparaissent dans les produits $\omega^{s+j,k-j}\wedge\omega^{s,j}_{\eta}$ on obtient les égalités suivantes:
\begin{eqnarray*}
(A_3)=g_{s+k}^{-1}\omega^{s+k+1,-1}\Big(\bar{\partial}_{_J } \eta^{s,k+1} +\sum_{j=0}^{k+1}\omega^{s+j,k+1-j}\wedge\eta^{s,j} -\sum_{j=-1}^{k} (-1)^{k+1-j}\eta^{s+j,k+1-j}\wedge \omega^{s,j}_{\eta}\Big)+ 
\\
\\
+g_{s+k}^{-1}\bar{\partial}_{_J }\omega^{s,k}+\sum_{j=-1}^{k} (-1)^{k-j}g_{s+k}^{-1}\omega^{s+j,k-j}\wedge\Big(\sum_{r=-1}^{j-1}(-1)^{j-r+1}\eta^{s+r,j-r}\wedge\omega^{s,r}_{\eta}+\omega^{s,j}\Big)\qquad
\\
\\
+\sum_{j=-1}^{k}\sum_{r=1}^{k-j} (-1)^{k-j}g_{s+k}^{-1}\omega^{s+j+r,k-j-r}\wedge\eta^{s+j,r}\wedge \omega^{s,j}_{\eta}=(A_4)\qquad\qquad\qquad\quad
\end{eqnarray*}
En faisant le changement d'indice $j'=j+r$, $r'=j$ dans la dernière somme double on a finalement:
\begin{eqnarray*}
(A_4)=g_{s+k}^{-1}\omega^{s+k+1,-1}\Big(\bar{\partial}_{_J } \eta^{s,k+1} +\sum_{j=0}^{k+1}\omega^{s+j,k+1-j}\wedge\eta^{s,j} -\sum_{j=-1}^{k} (-1)^{k+1-j}\eta^{s+j,k+1-j}\wedge \omega^{s,j}_{\eta}\Big)+ 
\\
\\
+g_{s+k}^{-1}\Big(\bar{\partial}_{_J }\omega^{s,k}+\sum_{j=-1}^{k} (-1)^{k-j}\omega^{s+j,k-j}\wedge\omega^{s,j}\Big)=
\omega^{s+k+1,-1}_{\eta}\wedge \omega^{s,k+1}_{\eta}\qquad\qquad\quad
\end{eqnarray*}
ce qui justifie la formule $(3^{s,k}_{\eta})$. On a alors qu'à la fin de cette récurrence toutes les matrices $\omega^{s,k}_{\eta} $ vérifient la relation $(3_{\eta})$.
Montrons maintenant que l'application $R$ est une action de semi-groupe . On se propose donc de montrer la formule $R(\eta_2,\omega_{\eta_1})=: \omega_{\eta_1,\eta_2}=\omega_{\eta_1\wedge\eta_2}$ qui en termes de composantes s'exprime sous la forme 
$\omega^{s,k}_{\eta_1,\eta_2}=\omega^{s,k}_{\eta_1\wedge\eta_2} $ pour tout $k\geq -1$. On montre la formule précédente par récurrence sur $k$. On remarque que la formule est évidente pour $k=-1$. En explicitant l'expression de $\omega^{s,k}_{\eta_1,\eta_2}$ et en utilisant l'hypothèse de récurence on a:
\begin{eqnarray*}
\omega^{s,k}_{\eta_1,\eta_2}= g_{s+k,2}^{-1}\Big(\bar{\partial}_{_J }\eta^{s,k}_2+\sum_{j=1}^{k+1}  \omega^{s+j,k-j}_{\eta_1}\wedge \eta^{s,j}_2-\sum_{j=-1}^{k-1} (-1)^{k-j}\eta^{s+j,k-j}_2\wedge \omega^{s,j}_{\eta_1\wedge \eta_2} +\omega^{s,k}_{\eta_1}\cdot g_{s,2}\Big) =
\\
\\
=(g_{s+k,1}\cdot g_{s+k,2})^{-1}\Big[g_{s+k,1}\cdot \bar{\partial}_{_J }\eta^{s,k}_2+\sum_{j=1}^{k+1}\bar{\partial}_{_J } \eta^{s+j,k-j}_1\wedge \eta^{s,j}_2 +\qquad\qquad\qquad\quad
\\
\\
\underbrace{+\sum_{j=1}^{k+1}\sum_{r=1}^{k-j+1}\omega^{s+j+r,k-j-r}\wedge\eta^{s+j,r}_1\wedge\eta^{s,j}_2+\sum_{j=1}^{k+1} \omega^{s+j,k-j}\wedge g_{s+j,1}\cdot\eta^{s,j}_2-}_{(1)} \qquad\qquad\quad
\\
\\
\underbrace{ -\sum_{j=1}^{k+1}\sum_{r=-1}^{k-j-1}(-1)^{k-j-r} \eta^{s+j+r,k-j-r}_1\wedge\omega^{s+j,r}_{\eta_1}\wedge\eta^{s,j}_2}_{(2)}
-\sum_{j=-1}^{k-1}(-1)^{k-j}g_{s+k,1}\cdot \eta^{s+j,k-j}_2\wedge\omega^{s,j}_{\eta_1\wedge\eta_2}+
\\
\\
+\Big(\bar{\partial}_{_J }\eta^{s,k}_1+\underbrace{ \sum_{j=1}^{k+1}  \omega^{s+j,k-j}\wedge \eta^{s,j}_1}_{(1)} -\underbrace{ \sum_{j=-1}^{k-1} (-1)^{k-j}\eta^{s+j,k-j}_1\wedge \omega^{s,j}_{\eta_1}}_{(2)}  +\omega^{s,k}\cdot g_{s,1}\Big)g_{s,2}\Big]  =(B_1)\quad 
\end{eqnarray*}
En rappelant l'expression du terme $\bar{\partial}_{_J}(\eta_1\wedge \eta_2)^{s,k} $, en faisant le changement d'indice $j'=j+r$, $r'=j$ dans la première  et deuxième somme double et en regroupant opportunément les termes on obtient:
\begin{eqnarray*}
(B_1)= (g_{s+k,1}\cdot g_{s+k,2})^{-1}\Big[\bar{\partial}_{_J}(\eta_1\wedge \eta_2)^{s,k}
\underbrace{ +\sum_{j=0}^{k+1}  \omega^{s+j,k-j}\wedge(\eta_1\wedge \eta_2)^{s,j}}_{(1)} +\omega^{s,k} -
\\
\\ 
-\sum_{j=-1}^{k-1}(-1)^{k-j}\eta^{s+j,k-j}_1\wedge\Big(\bar{\partial}_{_J } \eta^{s,j}_2    
\underbrace{+\sum_{r=0}^{j+1}\omega^{s+r,j-r}_{\eta_1} \wedge\eta^{s,r}_2+\omega ^{s,j}_{\eta_1}}_{(2)} \Big)-\qquad\quad
\\
\\
-\sum_{j=-1}^{k-1}(-1)^{k-j}g_{s+k,1}\cdot\eta^{s+j,k-j}_2\wedge\omega^{s,j}_{\eta_1\wedge\eta_2}\,\Big]=(B_2)\qquad\qquad\qquad\quad
\end{eqnarray*}
En utilisant l'hypothèse de récurrence et la définition des matrices $\omega^{s,j}_{\eta_1,\eta_2}$, on peut écrire le terme entre parenthèses rondes sous la forme suivante:
\begin{eqnarray*} 
\bar{\partial}_{_J } \eta^{s,j}_2  + \sum_{r=0}^{j+1}\omega^{s+r,j-r}_{\eta_1} \wedge\eta^{s,r}_2+\omega ^{s,j}_{\eta_1} =
\sum_{r=-1}^{j}(-1)^{j-r}\eta^{s+r,j-r}_2\wedge\omega^{s,r}_{\eta_1\wedge\eta_2}+\omega^{s,j}_{\eta_1\wedge\eta_2}
\end{eqnarray*}
On aura alors
\begin{eqnarray*} 
(B_2)=(g_{s+k,1}\cdot g_{s+k,2})^{-1}\Big[ \bar{\partial}_{_J}(\eta_1\wedge \eta_2)^{s,k}+\sum_{j=0}^{k+1}\omega^{s+j,k-j}\wedge(\eta_1\wedge \eta_2)^{s,j}
+\omega^{s,k} -\qquad\quad
\\
\\
-\sum_{j=-1}^{k-1}\sum_{r=-1}^{j}(-1)^{k-r}\eta^{s+j,k-j}_1\wedge\eta^{s+r,j-r}_2\wedge\omega^{s,r}_{\eta_1\wedge\eta_2}-\qquad\qquad\qquad\qquad
\\
\\
-\sum_{j=-1}^{k-1}(-1)^{k-j}\Big(\eta^{s+j,k-j}_1+ \eta^{s+j,k-j}_2\Big)\wedge\omega^{s,j}_{\eta_1\wedge\eta_2}-
\sum_{j=-1}^{k-1}(-1)^{k-j}\eta^{s+k,0}\wedge\eta^{s+j,k-j}_2\wedge\omega^{s,j}_{\eta_1\wedge\eta_2}\,\Big]
\end{eqnarray*}
En faisant le changement d'indice $j'=r$, $r'=j-r$ dans la somme double et en regroupant ce terme avec les deux dernières sommes, on obtient le terme cherché $\omega^{s,k}_{\eta_1\wedge\eta_2}$.  La transitivité de l'action $R$  est complètement claire par les  calculs qui ont permis de prouver que l'action même est bien définie\hfill $\Box$
\\
On va considérer maintenant quelques formules utiles pour le procédé itératif de la convergence rapide qui sera exposé en détail dans la sub-section (3.4). On explique formellement les  étapes  du procédé itératif. On désigne par $\omega_0:=\omega \in\Omega (U,\varphi,\psi,\bar{\partial})$ le choix initial de $\omega$. Au 
$k$-ème pas du procédé itératif on suppose avoir obtenu l'élément $\omega_k \in\Omega (U,[\varphi,\psi],\bar{\partial})$  et avoir déterminé le paramètre $\eta_{k+1}\in {\cal P}(U)$ en fonction de $\omega _k$. On définit alors l'élément $\omega _{k+1}:=R(\eta_{k+1},\omega _k)\equiv \omega _{k,\,\eta_{k+1} }$. Si on pose par définition
$$
\eta(k):={\displaystyle \bigwedge_{1\leq j\leq k }^{\longrightarrow}}\eta_j 
$$
on aura la formule  $\omega_k=\omega_{\eta(k)}$, grâce au fait que la recalibration $R$ est une action de semi-groupe. Si on pose par définition 
$g(k):=g(\eta(k))\in \Gamma(U)$ on aura que les composantes de $g(k)$ sont définies par la formule
$$
g_s(k) :={\displaystyle \prod_{0\leq j\leq k }^{\longrightarrow}}g_{s,j} 
$$
pour $s=0,...,m$, où $g_j:=g(\eta_j)$. On écrit maintenant, à l'aide de cette dernière définition, les composantes 
$\eta(k)^{s,t},\;t\geq 1$ du paramètre $\eta(k)$ défini précédemment, sous une forme utile pour la convergence vers une solution d'un problème différentiel qu'on exposera dans la troisième partie.
\begin{lem}. Pour tout entier $k \geq 1$ on peut écrire les composantes $\eta(k)^{s,t},\;t\geq 1$ du paramètre  $\eta(k)$ sous la forme
\begin{eqnarray}
\eta(k)^{s,t} =\Big(\sum_{\tau \in\Delta _t}\;\sum_{J\in J_k(\rho (\tau) )}\;{\displaystyle \bigwedge_{1\leq r\leq \rho (\tau) }^{\longrightarrow}} 
\;g_{s+\sigma' (\tau,r)}  (j_r-1)\cdot \eta^{s+\sigma (\tau,r),\,\tau_{\rho (\tau)+1-r } }_{j_r} \cdot g_{s+\sigma (\tau,r)} (j_r)^{-1}\,\Big)\cdot g_s(k)  
\end{eqnarray}
où
\begin{eqnarray*}
\Delta _t:=\left\{\tau\in\N^t\,|\,\tau_j\not=0 \Rightarrow \tau_{j-1} \not=0\,,\,\sum_{j=1}^t\tau_j =t\right\}
\\
\\  
\rho(\tau):=\max\{j\,|\,\tau_j\not=0\}\qquad\qquad\quad
\\
\\
J_k(\rho (\tau)):=\{J\in\{1,...,k\}^{\rho(\tau)}\,|\,j_1<...<j_{\rho(\tau)}\} 
\\
\\
\sigma'(\tau,r):=\sum_{j=1}^{\rho (\tau)+1-r}\tau_j \quad\text{et}\quad\sigma(\tau,r):=\sum_{j=1}^{\rho (\tau)-r}\tau_j 
\end{eqnarray*}
\end{lem} 
Remarquons que $J_k(\rho (\tau))=\emptyset$ si $k<\rho (\tau)$.
\\
$Preuve$. On remarque que l'expression (5) est évidente dans le cas $k=1,2$. Il est immédiat de vérifier à l'aide d'une récurrence élémentaire la validité de l'expression (5) pour $t=1$ et $k\geq 1$ entier quelconque. Dans ce cas la formule (5) s'écrit sous la forme
\begin{eqnarray*}
\eta(k)^{s,1}=\Big(\sum_{j=1}^kg_{s+1}(j-1)\cdot \eta^{s,1}_j \cdot g_s(j)^{-1}\Big)\cdot g_s(k)      
\end{eqnarray*}
On montre maintenant la validité de l'expression (5) en général à l'aide du procédé récursif suivant. On suppose vraie la formule (5) pour les composantes 
$\eta(k)^{\bullet,j} $, $j=1,...,t+1$ , $k \geq 1$ et on prouve la formule pour la composante $\eta(k+1)^{\bullet ,t+1}$ En effet on a 
\begin{eqnarray*}
\eta(k+1)^{s,t+1}:=(\eta(k)\wedge \eta_{k+1})^{s,t+1}:=\qquad\qquad\qquad
\\
\\
=\eta(k)^{s,t+1}\cdot g_{s,k+1}+g_{s+t+1} (k)\cdot \eta_{k+1}^{s,t+1}+      
\sum_{j=1}^t\eta(k)^{s+j,t+1-j}\wedge \eta_{k+1}^{s,j}=
\\
\\
=\Big(\sum_{\tau \in\Delta _{t+1} }\;\sum_{J\in J_k(\rho (\tau))}...\Big)\cdot g_s(k+1)  +g_{s+t+1} (k) \cdot \eta_{k+1}^{s,t+1}+\qquad      
\\
\\
\sum_{j=1}^t \sum_{\tau \in\Delta _{t+1-j} }\;\sum_{J\in J_k(\rho (\tau) )}\;{\displaystyle \bigwedge_{1\leq r\leq \rho (\tau) }^{\longrightarrow}} \;(...)\wedge
 g_{s+j} (k)\cdot  \eta_{k+1}^{s,j}=\qquad\quad
\end{eqnarray*}
\begin{eqnarray*}
= \Big(\sum_{j=1}^{k+1} g_{s+t+1}(j-1)\cdot \eta^{s,t+1}_j \cdot g_s(j)^{-1}\Big)\cdot g_s(k+1) +\qquad\qquad\qquad
\\
\\
+\Big({\displaystyle \sum_{\scriptstyle \tau \in\Delta _{t+1} 
\atop
\scriptstyle \rho (\tau) \geq 2}\;\sum_{J\in J_k(\rho (\tau))}...}\Big)\cdot g_s(k+1) + 
\Big({\displaystyle \sum_{\scriptstyle \tau \in\Delta _{t+1} 
\atop
\scriptstyle \rho (\tau) \geq 2}\;\sum_{\scriptstyle J\in J_{k+1} (\rho (\tau))
\atop
\scriptstyle  j_{\rho (\tau)}=k+1 }...}\Big)\cdot g_s(k+1) =\qquad
\\
\\ 
\Big(\sum_{j=1}^{k+1} g_{s+t+1}(j-1)\cdot \eta^{s,t+1}_j \cdot g_s(j)^{-1}\Big)\cdot g_s(k+1) + \Big({\displaystyle \sum_{\scriptstyle \tau \in\Delta _{t+1} 
\atop
\scriptstyle \rho (\tau) \geq 2}\;\sum_{J\in J_{k+1} (\rho (\tau))}...}\Big)\cdot g_s(k+1) 
\end{eqnarray*}
ce qui prouve la formule (5) pour la composante $\eta(k+1)^{s,t+1}$.\hfill $\Box$

%%%%%%%%%%%%%%%%%%%%%%%%%%%%%%%%%%%%%%%%%%%%%%%%%%%%%%%%%%%%%%%%%%%
\subsection{Troisième étape: formulation du problème différentiel.}
%%%%%%%%%%%%%%%%%%%%%%%%%%%%%%%%%%%%%%%%%%%%%%%%%%%%%%%%%%%%%%%%%%%

La partie principale de la preuve consiste à prouver l'existence, pour tout $x\in U$, d'un voisinage ouvert $V \subset U$ de $x$ et $g\in\Gamma(V)$ solution du système différentiel
$$
(\Sigma )\;
\left  \{
\begin{array}{lr}
\bar{\partial}(\psi\cdot g_0)=0 
\\
\\
\bar{\partial}_{_J}(g_{s-1}^{-1} \cdot \varphi _s \cdot g_s)=0  
\\
\\
s=1,...,m \qquad\qquad\qquad
\end{array}
\right.
$$
Bien évidement résoudre ce système différentiel équivaut à trouver une autre ${\cal E}$-résolution de  ${\cal G}_{_{|_V}}$ dans la classe $[\varphi ,\psi]$, à partir de la  ${\cal E}$-résolution donnée $(\varphi ,\psi)$, de telle sorte qu'elle admet des matrices de connexion $\omega^{s,0}$  nulles. 
Maintenant on va prouver les deux résultats suivants.
\begin{lem}
Pour tout choix de $\omega \in \Omega (U,\varphi ,\psi,\bar{\partial})$, l'existence d'une solution $g\in \Gamma(U)$ du système différentiel $(\Sigma)$ est équivalente à l'existence d'une solution $\eta\in {\cal P}(U),\; g=g(\eta)$, du système différentiel quasi-lineaire
$$
(S)\;
\left  \{
\begin{array}{lr}
\displaystyle{\bar{\partial}_{_J }\eta^{s,k}+\sum_{j=0}^{k+1}  \omega^{s+j,k-j}\wedge \eta^{s,j}+(-1)^{k}\eta^{s-1,k+1}\wedge \omega^{s,-1}_{\eta} +\omega^{s,k}=0}
\\
\\
k=0,...,m 
\\
s=0,...,m-k 
\end{array}
\right.
$$ 
\end{lem}  
La proposition suivante permet de traduire notre problème en termes purement différentiels.
\begin{prop}
Supposons données des matrices $\omega ^{s,k}\in M_{p_{s+k},p_s }({\cal E}_X^{0,k+1}(U))$, $s=0,...,m,\;k=-1,...,m-s,\;(s,k)\not=(0,-1)$ telles que $\omega ^{s-1,-1}\cdot\omega ^{s,-1}=0,\;s=2,...,m$ et $\bar{\partial}_{_J}\omega ^{s,-1}+\omega ^{s-1,0}\cdot \omega ^{s,-1}= \omega ^{s,-1}\cdot \omega ^{s,0},\;s=1,...,m$. Alors, pour $k\geq 0$, les relations $(3_{s,k})$
\begin{eqnarray*}
\bar{\partial}_{_J }\omega^{s,k} +\sum_{j=-1}^{k+1} (-1)^{k-j}  \omega^{s+j,k-j}\wedge \omega^{s,j}=0 
\end{eqnarray*}
constituent les conditions d'intégrabilité du système différentiel $(S)$.
\end{prop} 
On remarque que si  $\varphi_1 =0$ alors l'hypothèse d'intégrabilité $\bar{\partial}^{2}=0 $ s'exprime localement par la relation $\bar{\partial}_{_J} \omega ^{0,0} +\omega ^{0,0} \wedge \omega ^{0,0} =0$. Un résultat du à A.Grothendieck (voir \cite{ko-mal} que nous redemontrerons d'une façon différente), assure alors l'existence, pour chaque $x\in U$ d' un voisinage ouvert $V$ de $x$ et de  $h_0\in GL(p_0,{\cal E}(V))$ telle que $h^{-1}_0\bar{\partial}_{_J}h_0=\omega ^{0,0} $. Bien évidement cette équation est équivalente à l'unique équation qui constitue le système (S) dans le cas $m=0$, à condition d'identifier $h_0\equiv g_0(\eta)^{-1} $.
On obtient alors le diagramme commutatif suivant:
\begin{diagram}[height=1cm,width=1cm]
0&\rTo  &(Ker\,\bar{\partial})_{_{|_V}}&\rTo& {\cal G}_{_{|_V} }&\rTo^{\bar{\partial}} &{\cal G}_{_{|_V}}\otimes_{_{{\cal E}_{_V}}}{\cal E}^{0,1}_{_V}&
\\
&      &\uTo^{\psi_g}_{\wr}&            &\uTo^{\psi_g}_{\wr}&          &\uTo^{\wr}_{\psi_g \otimes \I_{(0,1)} }& 
\\
0&\rTo &{\cal O}^{\oplus p_0}_{_V}& \rTo &{\cal E}^{\oplus p_0}_{_V}& \rTo^{\bar{\partial}_{_J }} &({\cal E}^{0,1}_{_V})^{\oplus p_0}& 
\end{diagram} 
lequel permet de conclure dans ce cas. Venons-en maintenant à la preuve du lemme 3.3.1.
\\
\\
$Preuve\; du\; lemme\; 3.3.1$.
Soit $g\in \Gamma (U)$ une solution du système $(\Sigma)$. On écrit les composantes de $g$ sous la forme $g_s=\I_{p_s}+\eta^{s,0}$. Avec ces notations le système 
$(\Sigma)$ s'écrit sous la forme 
$$
(\Sigma_1)\;
\left  \{
\begin{array}{lr}
\bar{\partial}\psi_{\eta}=0 
\\
\\
\bar{\partial}_{_J }\omega ^{s,-1}_{\eta}=0  \qquad
\\
\\
s=1,...,m 
\end{array}
\right.
$$
On rappel que les calculs utilisés pour montrer que l'application de recalibration $R$ est bien définie, nous donnent les égalités
$$
\left  \{
\begin{array}{lr}
\bar{\partial}\psi_{\eta}=\psi \cdot ( \bar{\partial}_{_J }\eta^{0,0}+\omega ^{0,0}\wedge \eta^{0,0}+\omega ^{0,0})    
\\
\\
\bar{\partial}_{_J }\omega ^{s,-1}_{\eta}=-\omega ^{s-1,0}_{\eta}\cdot \omega ^{s,-1}_{\eta}+
\omega ^{s,-1}_{\eta}\cdot g_s^{-1} ( \bar{\partial}_{_J }\eta^{s,0}+\omega ^{s,0}\wedge \eta^{s,0}+\eta^{s-1,0}\wedge\omega ^{s,-1}_{\eta} +\omega ^{s,0})  
\\
\\
s=1,...,m 
\end{array}
\right.
$$
L'hypothèse d'exactitude de la ${\cal E}$-résolution locale implique alors l'existence d'une matrice $\eta^{0,1} \in M_{p_1,p_0}({\cal E}^{0,1}_X(U))$ telle que 
$$
\omega ^{0,0}_{\eta}= \bar{\partial}_{_J }\eta^{0,0}+\sum_{j=0}^1 \omega ^{j,-j}\wedge \eta^{0,j}+\omega ^{0,0}=0 
$$
(remarquons que la dépendance effective des matrices $\omega ^{\bullet,0}_{\eta}$ du paramètre $\eta$ est limitée aux composantes $\eta^{\bullet ,k},\;k=0,1$). On a donc que l'équation du système (S), relative aux indices $(s,k)=(0,0)$ est satisfaite. On obtient alors à l'aide d'une récurrence croissante sur les indices $s=1,...,m$ relatifs aux expressions précédentes des matrices $\bar{\partial}_{_J }\omega ^{s,-1}_{\eta}$ et de l'exactitude de la  ${\cal E}$-résolution locale, l'existence de matrices $\eta^{s,1} \in M_{p_{s+1} ,p_s}({\cal E}^{0,1}_X(U)),\;s=0,...,m$ telles que $\omega ^{s,0}_{\eta}=0$ pour les indices en question. Ces équations ne représentent rien d'autre que le système différentiel quasi-lineaire
$$
(S_1)\;
\left  \{
\begin{array}{lr}
\displaystyle{\bar{\partial}_{_J }\eta^{s,0}+\sum_{j=0}^{1}  \omega^{s+j,-j}\wedge \eta^{s,j}+\eta^{s-1,1}\wedge \omega^{s,-1}_{\eta} +\omega^{s,k}=0}
\\
\\
k=0,...,m 
\end{array}
\right.
$$ 
qui est évidement équivalent au système $(\Sigma)$. On rappelle aussi que les calculs relatifs à la bonne définition de l'application de recalibration $R$, nous donnent les égalités
\begin{eqnarray*}
\bar{\partial}_{_J }\omega^{s,k}_{\eta} +\sum_{j=-1}^{k} (-1)^{k-j}  \omega^{s+j,k-j}_{\eta}\wedge \omega^{s,j}_{\eta}=\qquad\qquad\qquad\qquad\qquad
\\
\\
=g^{-1}_{s+k} \omega^{s+k+1,-1}\Big(\bar{\partial}_{_J } \eta^{s,k+1} +\sum_{j=0}^{k+1}\omega^{s+j,k+1-j}\wedge\eta^{s,j} -\sum_{j=-1}^{k} (-1)^{k+1-j}\eta^{s+j,k+1-j}\wedge \omega^{s,j}_{\eta}+\omega ^{s,k+1} \Big)
\end{eqnarray*}
On obtient alors à l'aide d'une récurrence triangulaire, analogue à celle qui nous a permis de choisir les matrices $\omega ^{\bullet, \bullet}$ dans la première étape de la preuve, l'existence des matrices $\eta^{s,k}$ telles que $\omega^{s,k}_{\eta}=0$ pour $s=0,...,m,\;k=0,....,m$. Le système formé par ces équations est bien évidement équivalent au système (S), ce qui prouve le lemme.\hfill $\Box$
\\
Venons maintenant à la proposition 3.3. On commence par prouver la nécessité des conditions d'intégrabilité pour le système (S). La suffisance de ces conditions sera prouvée dans l'étape suivante de la preuve du théorème 1.
\\
$Preuve\;de \;la \;n\acute{e}cessit\acute{e}\; des \; conditions \;d'int\acute{e}grabilit\acute{e}\;pour \;le\; syst\grave{e} me\;(S).$
On commence par prouver la validité des relations $(3_{\bullet,k})$, $k\geq 0$ à l'aide d'une récurrence  croissante sur $k=0,..,m$. On rappelle que, le fait que par hypothèse on dispose des relations $\omega ^{s-1,-1}\cdot\omega ^{s,-1}=0,\;s=2,...,m$ et $(3_{\bullet,-1})$, combiné avec le fait que les équations du système (S), relatives aux indices $(s,0)$ ne raprésentent rien d'autre que les équations $\omega^{s,0}_{\eta}=0$, implique la validité des équations 
$\bar{\partial}_{_J }\omega^{s,-1}_{\eta}=0$. On suppose par hypothèse de récurence la validité des relations $(3_{\bullet,j})$ pour $j=-1,...,k-1$. En utilisant la validité des équations du système (S) et en appliquant l'opérateur $\bar{\partial}_{_J }$ à l'équation relative aux indices $(s,k)$ on obtient les égalités suivantes.
\begin{eqnarray*}
0=\bar{\partial}_{_J }\omega ^{s,k}\cdot g_s-(-1)^k \omega ^{s,k}\wedge \bar{\partial}_{_J }\eta ^{s,0} +
\sum_{j=1}^{k+1} \bar{\partial}_{_J } \omega ^{s+j,k-j}\wedge\eta^{s,j}-\quad
\\
\\
-\sum_{j=1}^{k+1} (-1)^{k-j} \omega ^{s+j,k-j}\wedge \bar{\partial}_{_J } \eta^{s,j} +(-1)^k\bar{\partial}_{_J }\eta^{s-1,k+1}\wedge \omega^{s,-1}_{\eta}=\quad\;
\\
\\
=\bar{\partial}_{_J }\omega ^{s,k}\cdot g_s-\sum_{j=1}^{k+1}\;\sum_{r=-1}^{k-j+1}(-1)^{k-j-r}\omega ^{s+j+r,k-j-r}\wedge \omega ^{s+j,r}\wedge \eta^{s,j}-
\\
\\
-\sum_{j=1}^{k+1} (-1)^{k-j} \omega ^{s+j,k-j}\wedge \bar{\partial}_{_J } \eta^{s,j} +(-1)^k\bar{\partial}_{_J }\eta^{s-1,k+1}\wedge \omega^{s,-1}_{\eta} \qquad 
\end{eqnarray*}
En faisant le changement d'indice $j'=j+r$, $r'=j$ dans la somme double on obtient
\begin{eqnarray*}
0=\bar{\partial}_{_J }\omega ^{s,k}\cdot g_s-\sum_{j=0}^{k+1}(-1)^{k-j}
 \omega ^{s+j,k-j}\wedge\Big(\bar{\partial}_{_J } \eta^{s,j}+\sum_{r=1}^{j+1}\omega ^{s+r,j-r}\wedge \eta^{s,r}\,\Big)+
\\
\\
+(-1)^k\bar{\partial}_{_J }\eta^{s-1,k+1}\wedge \omega^{s,-1}_{\eta} =\Big(\bar{\partial}_{_J }\omega^{s,k} +\sum_{j=0}^{k+1} (-1)^{k-j}  \omega^{s+j,k-j}\wedge \omega^{s,j}\Big)\cdot g_s+
\\
\\
+(-1)^k\Big(\bar{\partial}_{_J }\eta^{s-1,k+1}+\sum_{j=0}^{k+1}\omega ^{s+j,k-j}\wedge \eta^{s-1,j+1}\Big)\wedge \omega^{s,-1}_{\eta}=\qquad\quad
\end{eqnarray*}
\begin{eqnarray*}
=\Big(\bar{\partial}_{_J }\omega^{s,k} +\sum_{j=-1}^{k+1} (-1)^{k-j}  \omega^{s+j,k-j}\wedge \omega^{s,j}\Big)\cdot g_s 
\end{eqnarray*}
L'inversibilité de $g_s$ permet alors de conclure la  preuve de la nécessité des conditions d'intégrabilité $(3_{\bullet,k})$, $k\geq 0$ pour le système $(S)$. \hfill $\Box$
\\
La proposition (3.3) nous suggère de considérer les définitions suivantes. On définit l'ensemble 
\begin{eqnarray*} 
\Omega (U,p)\subset {\displaystyle  \bigoplus_
{{\scriptstyle  s=0,...,m
\atop
\scriptstyle  k=-1,...,m-s}
\atop
\scriptstyle (s,k)\not=(0,-1)} }M_{p_{s+k} ,p_s}({\cal E}^{0,k+1}_X(U))
\end{eqnarray*} 
$p=(p_1,...,p_m)$, constitué des éléments $\omega =(\omega ^{s,t})_{s,t}$ tels que $\omega^{s-1,-1}\cdot \omega ^{s,-1}=0$ et la relation (3) soit satisfaite. Si on dispose de matrices $\omega^{s,-1}_0\in M_{p_{s-1},p_s}({\cal E}_X(U)) $, $s=1,...,m$ qui vérifient la relation écrite précédemment on peut définir l'ensemble
$$
\Omega (U,\omega ^{\bullet,-1}_0):=\{\omega \in \Omega (U,p)\,|\,\exists g\in \Gamma(U):\omega ^{\bullet,-1}=\omega ^{\bullet,-1}_{0,\,g}\}  
$$
Les calculs relatifs à la proposition (3.2) permettent d'étendre la recalibration $R$ à l'application
$$
R:{\cal P}(U)\times \Omega (U,\omega ^{\bullet,-1}_0)     \longrightarrow \Omega (U,\omega ^{\bullet,-1}_0)
$$
laquelle est encore une action de semi-groupe. A partir de maintenant on va considérér plus généralement la recalibration en termes de l'application définie précédemment.
Venons-en maintenant  à un  préliminaire technique avant d'exposer la preuve de l'existence des solutions pour le système différentiel $(S)$.

%%%%%%%%%%%%%%%%%%%%%%%%%%%%%%%%%%%%%%%%%%%%%%%%%%%%%%%%%%%%%%
\subsection*{Choix des normes et opérateur de Leray-Koppelman.}
%%%%%%%%%%%%%%%%%%%%%%%%%%%%%%%%%%%%%%%%%%%%%%%%%%%%%%%%%%%%%% 

A partir de maintenant on va supposer que $U=B_1(0)$ est la boule ouverte de $\C^n$ de centre l'origine et de rayon unité. 
Si $A\in M_{k,l}(\C)$  on définit la norme $\|A\|:=\sup_{v\in\C^l-\{0 \} }\|Av\|/\|v\|$ et si $u=\sum'_{|I|=q}u_I\, d\bar{z}_I$ est une $(0,q)$-forme à coefficients des $(k,l)$-matrices à coefficients dérivables jusque à l'ordre $h\geq 0$, on définit une norme de Hölder invariante par changement d'échelle
$$
\|u\|_{r,\,h,\,\mu ,\,q}:=\sum _{\displaystyle
\scriptstyle   |I|=q 
\atop
\scriptstyle |\alpha |\leq h }
S_{|\alpha |} \,r^{|\alpha |+q}\,\|\partial^\alpha u_I\|_{r,\mu } 
$$    
où 
$$
\|f\|_{r,\,\mu}:=\sup_{z\in B_r} \,\|f(z)\|+\sup_{\displaystyle
\scriptstyle  z,\zeta \in B_r
\atop
\scriptstyle  z\neq \zeta }
 r^{\mu}\,\frac{ \|f(z) -f(\zeta ) \|}{\|z-\zeta \|^{\mu } } 
$$    
avec $\mu \in (0,1)$ une constante fixé une fois pour toutes dans notre problème et $(S_k)_{k\geq 0}\subset (0,\infty),\;S_0:=1$ est une suite de réels qu'on construira ensuite et qui vérifie l'inégalité
$$
S_k\leq\left[\max_{|\alpha +\beta |=k}{\alpha +\beta\choose \alpha  }\right]^{-1} S_j\,S_{k-j}   
$$
pour tout $k\geq 1$ et $j=1,...,k-1$ On remarque que si le dégrée $q\geq 1$ on a que la norme $\|u\|_{r,h,q}$ tend vers zéro lorsque le rayon $r$ tend vers zéro. On désignera par ${\cal C}^{h,\mu} _{0,q}(\bar{B}_r, M_{k,l}(\C))$ l'espace de Banach de $(0,q)$-formes sur la boule fermé $\bar{B}_r $, à coefficients des $(k,l)$-matrices à coefficients dérivables jusque à l'ordre $h\geq 0$, telles que la norme $\|\cdot\|_{r,h,q}$ soit finie (on ne notera pas les dimensions des matrices).  On remarque que si $u\in {\cal C}^{h,\mu } _{0,q}(\bar{B}_r , M_{k,l}(\C))$ et  $v\in {\cal C}^{h,\mu } _{0,p}(\bar{B}_r, M_{l,t}(\C))$ alors on a l'inégalité $\|u \wedge v\|_{r,h ,p+q}\leq \|u\|_{r,h,q}\cdot\|v\|_{r,h ,p}$. On rappelle  très rapidement la définition de l'opérateur de Leray-Koppelman (voir les ouvrages classiques de Henkin-Leiterer \cite{He-Le}, de Range \cite{Ra} et l'article de Harvey-Polkin \cite{Ha-Po}). L'opérateur de Leray-Koppelman de la boule unité 
$$
T_q: {\cal C}^{h,\mu } _{0,q+1}(\bar{B}_1, M_{k,l}(\C))\longrightarrow  {\cal C}^{h,\mu } _{0,q}(B_1, M_{k,l}(\C))
$$
pour $q\geq 0$ est défini par une formule du type
$$
T_q\,u(z):=\int\limits_{\zeta \in B_1} u(\zeta )\wedge K_q (\zeta ,z)+ 
\int\limits_{\zeta \in \partial B_1} u(\zeta )\wedge k_q (\zeta ,z)
$$
où le premier opérateur intégral s'exprime en termes des coefficients $u_I$ de la forme $u$, par des termes du type
$$
\int\limits_{\zeta \in B_1} u_{I} (\zeta )\cdot K (\zeta ,z)\,d\lambda (\zeta )\qquad\text{avec}\qquad K (\zeta ,z)=\frac{\bar{\zeta }_l-\bar{z}_l  }{|\zeta -z|^{2n} }  
$$
et le deuxième par des termes du type
$$
\int\limits_{\zeta \in \partial B_1} u_{I} (\zeta )\cdot k (\zeta ,z)\,d\sigma  (\zeta )\qquad\text{avec}\qquad k(\zeta ,z)=
\frac{(\bar{\zeta }_j-\bar{z}_j)\cdot\bar{\zeta }_k  }{|\zeta -z|^{2l+2}\cdot\left[\bar{\zeta }\cdot(\bar{\zeta }-\bar{z})\right]^{n-1-l}  }  
$$
$ l=0,...,n-2$. L'opérateur de Leray-Koppelman $T_{r,q} : {\cal C}^{h,\mu} _{0,q+1}(\bar{B} _r, M_{k,l}(\C))\longrightarrow  {\cal C}^{h,\mu } _{0,q}(B_r, M_{k,l}(\C))$, $q\geq 0$ de la boule de rayon $r$ et de centre l'origine est défini par la formule $T^{r,q}:=(\lambda _r^{-1})^*\circ T_q \circ \lambda ^*_r$ avec $\lambda _r:B_1\rightarrow B_r$ l'homothétie de rapport $r$. Les propriétés de l'opérateur de Leray-Koppelman qui nous intéressent sont les suivantes:
\\
1) Pour toute forme différentielle $u\in {\cal C}^{h,\mu } _{0,q+1}(\bar{B} _r , M_{k,l}(\C))$ on a la formule d'homotopie:
\begin{eqnarray}
u= \bar{\partial}_{_J }T_{r,q}\, u+T_{r,q+1}\,\bar{\partial}_{_J }u
\end{eqnarray}
valable sur la boule $B_r$. 
\\
2) Il existe une suite de poids $S=(S_k)_{k\geq 0}$ de la norme de Hölder introduite précédemment telle que pour toute forme différentielle $u\in {\cal C}^{h,\mu } _{0,q+1}(\bar{B} _r , M_{k,l}(\C))$ on a l'estimation intérieure:
\begin{eqnarray}
\|T_{r,q} \,u\|_{r(1-\sigma ),\,h+1,\,\mu ,\,q} \leq C\cdot \sigma ^{-s(h)}\cdot\|u\|_{r,\,h,\,\mu ,\,q+1 }   
\end{eqnarray}
avec $\sigma \in (0,1)$, $s(h)=2n+k+2$ et $C=C(n,\mu )>0$ une constante $ind\acute{e} pendente$ de la régularité $h$. La preuve de la formule d'homotopie est exposé dans les ouvrages classiques mentionnés précédemment. Une estimation analogue à la (7) à  déjà été montrée par S.Webster (voir \cite{We-1}). Nous utiliserons essentiellement les mêmes arguments de Webster pour montrer celle ci. La différence avec l'estimation obtenue par Webster consiste dans le fait que la constante $C>0$ est indépendante de la régularité $h$. A partir de maintenant on désignera par $C$ une constante strictement positive indépendante de la régularité des formes. Pour prouver l'estimation (7) il suffit de se restreindre au cas $r=1$, la norme étant choisie invariante pour changement d'échelle. On considère à ce propos une fonction $\rho \in  {\cal C}^{\infty}(\R,\,[0,1])$ telle que $\rho (x)=1$ pour $x\leq 0$ et $\rho (x)=0$ pour $x\geq 1$. On définit alors la fonction de cutt-off $\chi_{\sigma }$, avec $\sigma \in(0,1)$, par la loi 
$$
\chi _{\sigma } (z):=\left  \{
\begin{array}{lr}
1  &  \text{si}\quad  |z|\leq 1-\sigma/2 
\\
\rho (2\sigma ^{-1}(|z|-1+\sigma/2 )) & \text{si} \quad 1-\sigma/2\leq |z|
\end{array}
\right.
$$ 
On aura alors, comme conséquence de l'invariance par translation de $ K(\zeta ,z)$, les égalités suivantes: 
\begin{eqnarray*}
J_1^{\alpha }(z):= \partial^{\alpha }_z \int\limits_{\zeta \in B_1} u_I (\zeta )\cdot K (\zeta ,z)\,d\lambda (\zeta )=\qquad\qquad\qquad\qquad\qquad
\\
\\
=\partial^{1_{\alpha } } _z \int\limits_{\zeta \in B_1} \partial^{\alpha-1_{\alpha }  }_{\zeta }(\chi_{\sigma }\cdot u_I) (\zeta )\cdot K (\zeta ,z)\,d\lambda (\zeta )+
 \int\limits_{\zeta \in B_1-B_{1-\sigma/2 } }((1-\chi_{\sigma })\cdot u_I) (\zeta )\cdot\partial^{\alpha }_z K (\zeta ,z)\,d\lambda (\zeta )
\end{eqnarray*}
pour $|z|\leq 1-\sigma$ et pour tout multi-indice $|\alpha |=k+1, \;k\geq 0$. Ici on désigne par $1_{\alpha}$ un multi-indice tel que $|1_{\alpha}|=1$ et $1_{\alpha }\leq \alpha $. La théorie classique du potentiel (voir par exemple \cite{Gi-Tru}) nous fournit alors les estimations 
\begin{eqnarray*}
\left\| \partial^{1_{\alpha } } _z \int\limits_{\zeta \in B_1} \partial^{\alpha-1_{\alpha }}_{\zeta }(\chi_{\sigma }\cdot u_I) (\zeta )\cdot K (\zeta ,z)\,d\lambda (\zeta )\right\|_{1-\sigma ,\,\mu } \leq C(n,\mu )\cdot\sigma^{-2n-1}\cdot\|\partial^{\alpha-1_{\alpha }}(\chi_{\sigma }\cdot u_I)\|_{1 ,\,\mu }  
\end{eqnarray*}
et
$
|\partial^{\alpha }_z K (\zeta ,z)|\leq  C(n,\mu )\cdot|\zeta -z|^{1-2n-|\alpha |}  
$
. On aura alors l'estimation  suivante:
\begin{eqnarray*}
\|J_1^{\alpha } \|_{1-\sigma ,\,\mu } \leq C(n,\mu)\cdot \sigma ^{-2n-1}  \cdot \|\partial^{\alpha-1_{\alpha}}(\chi_{\sigma }\cdot u_I)\|_{1,\,\mu }   
+C(n,|\alpha |)\cdot(\sigma/2)^{1-2n-|\alpha |} \cdot\|u_I \|_{1,\,0} 
\end{eqnarray*}
et donc
\begin{eqnarray*}
\|J_1\|_{1-\sigma ,\,h+1,\,\mu } \leq \sum_{|\alpha |\leq h+1}S_{|\alpha |}  \|J_1^{\alpha } \|_{1-\sigma, \,\mu }\leq \qquad\qquad\qquad\qquad
\\
\\
\leq\|J_1\|_{1-\sigma ,\,\mu } +  C(n,\mu )\cdot \sigma ^{-2n-1} \cdot\sum_{|\alpha |\leq h} \sum_{\beta \leq\alpha }S_{|\alpha |}{\alpha \choose\beta }  
\,\sigma ^{-|\beta |}\|\partial^{\beta  }\rho  \|_{1,\,\mu } \cdot\|\partial^{\alpha-\beta  }u_I \|_{1,\,\mu }+
\\
\\
+\sigma ^{-2n-h} \cdot\|u_I \|_{1,\,0}\cdot\sum_{|\alpha |\leq h+1} S_{|\alpha |} \cdot C(n,|\alpha |)\cdot 2^{|\alpha |+2n-1} \qquad\qquad\qquad
\end{eqnarray*}
Pour un choix convenable de la suite $ S:=(S_k)_{k\geq 0}$, qu'on présentera ensuite, on peut se ramener à supposer que $\|\rho  \|_{1,\,S,\,\mu }:=\sum_{\alpha \geq 0}S_{|\alpha |}\|\partial^{\alpha  }\rho  \|_{1,\,\mu }<+\infty $ et $\sum_{\alpha \geq 0}S_{|\alpha |} \cdot C(n,|\alpha |)\cdot 2^{|\alpha |+2n-1} <+\infty  $. On aura alors l'estimation 
\begin{eqnarray}
\|J_1^{\alpha } \|_{1-\sigma ,\,h+1,\,\mu }
\leq C(n,\mu ) \cdot\sigma ^{-2n-1}\cdot\|u_I \|_{1,\,\mu } +C(n,\mu ) \cdot\sigma ^{-2n-h-1}\cdot \|\rho \|_{1,\,h,\,\mu } \cdot\|u_I \|_{1,\,h,\,\mu }+\nonumber
\\\nonumber
\\
+C\cdot \sigma ^{-2n-h} \cdot\|u_I \|_{1,\,0} \leq C \cdot\sigma ^{-2n-h-1}\cdot\|u_I \|_{1,\,h,\,\mu}\qquad\qquad\qquad\qquad
\end{eqnarray}
Enfin pour estimer les termes du type
\begin{eqnarray*}
J_2^{\alpha }(z):= \partial^{\alpha }_z \int\limits_{\zeta \in\partial  B_1} u_I (\zeta )\cdot k (\zeta ,z)\,d\sigma  (\zeta )=\int\limits_{\zeta \in\partial  B_1} u_I (\zeta )\cdot\partial^{\alpha }_z k (\zeta ,z)\,d\sigma  (\zeta )
\end{eqnarray*}
 avec $|z|\leq 1-\sigma $, il suffit de dériver $|\alpha |+1$ fois le noyau $k(\zeta ,z)$, de remarquer l' estimation élémentaire:
\begin{eqnarray*}
\frac{|J_2^{\alpha }(z)-J_2^{\alpha }(\bar{z})| }{|z-\bar{z}|^{\mu } } \leq |z-\bar{z}|^{1-\mu }\cdot
\int\limits_{\zeta \in\partial  B_1} |u_I (\zeta )|\cdot|\nabla_z \,\partial^{\alpha }_z k (\zeta ,\hat{z}_{\zeta}  )|  \,d\sigma  (\zeta )
\end{eqnarray*}
(où $\hat{z}_{\zeta}$ est un point entre $z$ et $\bar{z} $) et les inégalités $|\zeta -z|\geq 3\sigma $, $ |\bar{\zeta }\cdot(\zeta -z)|\geq 3\sigma  $ pour $|\zeta |=1$ et $|z|\leq 1-3\sigma $. 
On aura alors l'estimation  
$$
\|J_2^{\alpha } \|_{1-\sigma ,\,\mu } \leq C(n,|\alpha |) \cdot \sigma^{-2n-1-|\alpha |}\cdot\|u_I \|_{1,0}
$$
Par l'hypothèse faite précédemment sur la suite de poids on aura l'estimation 
$$
\|J_2^{\alpha } \|_{1-\sigma ,\,h+1,\,\mu } \leq C\cdot \sigma ^{-s(h)}\cdot\|u_I \|_{1,0}
$$
laquelle combinée avec l'estimation (8) nous donne l'estimation (7) sur la boule de rayon unité. 
\\
Venons-en à la définition de la suite $S$ laquelle sera déterminée en partie par l'exigence de satisfaire les hypothèses faites dans les calculs précédents. On pose par définition
$$
A_k:=\sum_{|\alpha |=k}\max\{\|\partial^{\alpha}\rho \|_{1,0},\,\|\partial^{\alpha}\omega ^{s,-1}  \|_{1,0},\,s=0,..,m,\}  
$$
$B_k:=(\max\{A_k,C(n,k)\})^{-1}$ si $\max\{A_k,C(n,k)\}\not=0$  et 1 sinon, $D_k:=[ \max_{|\alpha +\beta |=k}{\alpha +\beta \choose\alpha  }]^{-1}$. On pose par définition $S_0=1,\;S_1=B_1>0$ et on définit $S_k$, $k\geq 2$ à l'aide de la formule récursive
$$
0<S_k:=\min\{2^{-k}B_k,\,R_k,\,L_k,\,D_k\cdot\min_{1\leq j\leq k-1}S_j\cdot S_{k-j}\} 
$$ 
où $R_k$, $L_k$ sont des constantes qui seront déterminées dans la deuxième étape. On désigne par $S(\omega )$ la suite de poids obtenue si on pose 
$R_k=L_k=+\infty$, dans la définition précédente des poids.
Avec ces définitions on aura $\|\omega ^{\bullet,-1}\|_{1,S}\leq \|\omega ^{\bullet,-1}\|_{1,S(\omega )}<+\infty$. Pour simplifier les notations on identifiera 
$\|\cdot\|_{r,h,\mu ,q}\equiv \|\cdot\|_{r,h}$  et $\|\partial^h f\|_{\bullet}\equiv\sum_{|\alpha |=h}\|\partial^{\alpha }  f\|_{\bullet}$ dans la suite. 
\\
\\
%%%%%%%%%%%%%%%%%%%%%%%%%%%%%%%%%%%%%%%%%%%%%%%%%%%%%%%%%%%%%%%%%%%%%%%%%%%%%%%%%%%%%%%%%%%%%%%%%%%%%%%%%%%%%%%%%%%%%%%%%%%%%%%%%%%%%%%%%%%%%%%%%%%%%%%%%%%
\subsection{Quatrième étape: présentation du schéma de convergence rapide de type Nash-Moser et existence d'une solution du problème différentiel (S).}

\subsubsection{Estimation fondamentale du schéma de convergence rapide. } 
%%%%%%%%%%%%%%%%%%%%%%%%%%%%%%%%%%%%%%%%%%%%%%%%%%%%%%%%%%%%%%%%%%%%%%%%%%%%%%%%%%%%%%%%%%%%%%%%%%%%%%%%%%%%%%%%%%%%%%%%%%%%%%%%%%%%%%%%%%%%%%%%%%%%%%%%%%%

Dans cette partie de la preuve on va montrer l'existence d'un paramètre de recalibration $\eta$ des éléments $\omega$, lequel permettra un contrôle quadratique de la norme des matrices $\omega ^{\bullet ,t}_{\eta},\;t\geq 0$ en termes de la norme des matrices $\omega ^{\bullet ,t},\;t\geq 0$. Ce contrôle est essentiel pour montrer la convergence vers zéro de la norme des matrices $\omega ^{\bullet ,t}_k,\;t\geq 0$ obtenues au $k$-ème pas du procédé itératif de la convergence rapide. La convergence vers une solution du problème différentiel (S) est appelé rapide en raison de de l'estimation quadratique mentionné précédemment. Avant de prouver l'estimation en question on va introduire quelques notations utiles pour la suite. Soit $\omega \in \Omega (B_1, p)$. Pour $r\in (0,1)$ on définit les quantités  
\begin{eqnarray*} 
a_h(\omega ,r):=\max\{\|\omega ^{s,k}_{\eta}\|_{r,h}\,|\,0\leq s \leq m\,,0\leq k \leq m-s\}
\\
\\
c(\omega ):=\max\{\|\omega ^{s,-1}\|_{1,\,S(\omega )}\,|\,1\leq s \leq m\}\qquad\quad
\end{eqnarray*} 
On remarque que, par définition de la norme de Hölder, la quantité $a_h(\omega ,r)$ tend vers zéro lorsque le rayon $r$ tend vers zéro. Pour tout $\sigma \in (0,1)$ on définit les rayons
$
r_l:=r(1-l\cdot \sigma _m)
$
pour $l=0,...,m+1$ où on pose par définition $\sigma _m:=\sigma /(m+1)$. Ensuite on définit par récurrence décroissante sur $k=m,...,0$, les constantes $L_k=L_k(C,c(\omega ))>0$ par les formules $L_m:=C$ et $L_{k-1}:=\max\{C, 2c(\omega )\cdot C\cdot L_k\}$. A partir de maintenant on désigne par $\varepsilon \in (0,1/2)$ une constante fixée telle que pour toutes les matrices $A\in M_{p_s,p_s}(\C)$ telle que $\|A\|<\varepsilon$ on a l'inversibilité de la matrice $\I_{p_s}+A $. Avec ces notations on a la proposition suivante.
\begin{prop}
Supposons donné $\omega \in \Omega (B_1, p),\;\sigma \in (0,1),\;r\in (0,1),\;h\in \N$ et les poids $0<S_j\leq S_j(\omega ),\;j=0,...,h+1$ de la norme de Hölder 
$\|\cdot\|_{r,h+1}$. Supposons que le rayon $r$ soit suffisamment petit pour assurer les estimations 
$$
L_0(C,c(\omega))\cdot \sigma _m^{-(m+1)\cdot s(h)}\cdot a_h(\omega ,r)<\varepsilon  
$$
et $a_h(\omega ,r)\leq 1$, où la quantité $a_h(\omega ,r)$ est calculée par rapport au poids $S_j,\;j=0,...,h$. Supposons de plus que  le poids  $S_{h+1}$ soit suffisement petit pour pouvoir assurer l'estimation:
\begin{eqnarray*}
S_{h+1} \|\partial^{h+1} \omega^{s,k}_I \|_{r ,\mu } \leq \| \omega^{s,k}_I \|_{r ,\mu }
\end{eqnarray*}
pour tout $k=0,...,m,\;s=0,...,m-k$ et $|I|=k+1$. Si on définit les composantes  $\eta^{s,k}$ du paramètre de recalibration 
 $\eta\in {\cal P}(\bar{B} _{r(1-\sigma )})$ par la formule de récurrence décroissante sur 
$k=m,...,0$
$$
\eta^{s,k}:=-T_{r_{m-k},\,k} \,\Big(\omega ^{s,k}+\omega ^{s+k+1,-1}\wedge\eta^{s,k+1}+(-1)^k\eta^{s-1,k+1} \wedge\omega ^{s,-1}\Big)\in M_{p_{s+k},p_s}({\cal E}^{0,k}_X(B_{r_{m-k}} )  )     
$$ 
pour $s=0,...,m-k$, alors on aura la validité des  estimations
\begin{eqnarray}
\|\eta^{\bullet,k}  \|_{r(1-\sigma ) ,h+1 } \leq L\cdot\sigma_m^{-(m+1)\cdot s(h)}\cdot a_h(\omega ,r)   
\\\nonumber
\\
\|\omega^{\bullet,k}_{\eta} \|_{r(1-\sigma ),h+1} \leq R \cdot\sigma_m ^{\nu (m,h)}\cdot a_h(\omega ,r) ^2 \quad\;
\end{eqnarray}
 pour tout $k= 0,...,m$, avec $L=L(C,c(\omega )):=\max_k L_k\,>0,\;R=R(C,c(\omega ))>0$ constantes positives et $\nu (m,h):=[(m+2)\cdot m+1]\cdot s(h)$, $(m\geq 0,\,h\geq 0)$. 
\end{prop} 
$Preuve$.
Dans les calculs qui suivront on utilisera les identifications $a_h\equiv a_h(\omega ,r)$ et $c=c(\omega) $. On commence par prouver l'estimation (9). Si on définit $\sigma_{m,l}>0$ par la formule $r_{l+1}=r_l(1+\sigma_{m,l})$ on a que $\sigma_{m,l}\geq \sigma _m$. On obtient alors à l'aide de l'estimation (7) et  d'une récurrence élémentaire décroissante sur $k=m,...,0$, l'estimation suivante
\begin{eqnarray}
\|\eta^{\bullet,k}  \|_{r_{m-k+1} ,h+1 } \leq L_k\cdot\sigma_m^{-(m-k+1)\cdot s(h)}\cdot a_h   
\end{eqnarray}
laquelle prouve l'estimation (9). Venons maintenant à la preuve de l'estimation (10). Pour cela on définit les troncatures 
$\eta_{[t]}  :=({\eta_{[t]}}^{s,k} )_{s,k}\in{\cal P}(B_{r_{m-t}}) $ du paramètre $\eta$ définit dans l'hypothèse de la proposition (3.4), de la façon suivante;  
${\eta_{[t]}}^{s,k} :=0$ si $k<t$ et  
${\eta_{[t]}}^{s,k}  :=\eta^{s,k}$ sur la boule $B_{r_{m-t}}$, si $k\geq t$. Par définition de la recalibration avec paramètre $\eta_{[k+1]}$ on aura alors: 
\begin{eqnarray*}
\omega ^{s,k}_{\eta_{[k+1]}}=  \omega^{s,k}+\omega^{s+k+1,-1} \wedge \eta^{s,k+1}+(-1)^k\eta^{s-1,k+1} \wedge\omega ^{s,-1}
\end{eqnarray*}
sur la boule $B_{r_{m-k-1}} $ pour $k=0,...,m$. 
On montre maintenant  à l'aide d'une récurrence en ordre décroissant sur $k=m,...,0$, l'estimation quadratique 
\begin{eqnarray}
\|\omega^{\bullet,k}_{\eta_{[k]} } \|_{r_{m-k+1} ,h+1} \leq R_k\cdot\sigma_m^{-b(m,k,h)} \cdot a_h^2 
\end{eqnarray}
où $R_k=R_k(C,c)>0$ est une constante positive, $b(m,k,h):= 2(m-k+1)\cdot s(h)$ si $k\geq 1$ et $b(m,0,h):= 2m+1$. Le diagramme suivant montre les matrices qui interviennent dans la définition de la matrice $\omega^{1,2}_{\eta_{[2]}}$ dans le cas  $m=4$.
%%%%% ICI
\begin{figure}[hbtp]
\begin{center} 
\input dgt5.pstex_t
    \caption{}
    \label{fig5}
\end{center} 
\end{figure}    
%%%%%%%%
On commence par prouver l'estimation (12) pour les indices $k=m,...,1$, (on suppose donc $m\geq 1$ dans les calculs qui suivront). En utilisant la formule d'homotopie pour l'opérateur $\bar{\partial}_{_J }$ et la condition 
$(3_{0,m})$ on a, pour $k=m$, les égalités suivantes:
\begin{eqnarray*}
\omega^{0,m}_{\eta_{[m]} }= \bar{\partial}_{_J }\eta^{0,m}+ \omega^{m,0}\wedge \eta^{0,m}- (-1)^m\eta^{0,m}\wedge \omega^{0,0}_{\eta_{[m]} } +\omega^{0,m}=\qquad\qquad\quad
\\
\\
=T_{r,m+1}\,\bar{\partial}_{_J }\omega^{0,m}-\omega^{m,0}\wedge T_{r,m} \,\omega^{0,m}+(-1)^m(T_{r,m}\,\omega^{0,m})\wedge \omega^{0,0}_{\eta_{[m]} }=\qquad\qquad
\\
\\
= -\sum_{j=0}^{m} (-1)^{m-j} T_{r,m+1}\, (\omega^{j,m-j}\wedge \omega^{0,j})
 -\omega^{m,0}\wedge T_{r,m} \,\omega^{0,m}+(-1)^m(T_{r,m}\,\omega^{0,m})\wedge \omega^{0,0}_{\eta_{[m]} }
\end{eqnarray*}
Le diagramme suivant montre les matrices qui interviennent dans la relation (3) dans le cas où $k=m$.
%%%%% ICI
\begin{figure}[hbtp]
\begin{center} 
\input dgt6.pstex_t
    \caption{}
    \label{fig6}
\end{center} 
\end{figure}    
%%%%%%%%
\\ 
On remarque que la relation $(3_{0,m})$ est la seule, parmi les autres relations $(3_{\bullet,\bullet})$, qui ne présente pas de facteurs de type $\omega^{\bullet,-1} $ dans les termes quadratiques.         
On estime maintenant la norme de la matrice $\omega^{0,m}_{\eta_{[m]} }$ à l'aide de l'expression précédente. On obtient l'inégalité suivante:
\begin{eqnarray*}
\|\omega^{0,m}_{\eta_{[m]} }\|_{r_1, h+1} \leq C\cdot \sigma_m ^{-s(h)}\cdot a_h^2 +\|\omega^{m,0}\wedge T_{r,m} \,\omega^{0,m} \|_{r_1, h+1}+
\\
\\
+\| (T_{r,m}\,\omega^{0,m})\wedge \omega^{0,0}\|_{r_1, h+1}  +c\cdot\| \eta^{0,m} \|^2_{r_1, h+1}\qquad\quad
\end{eqnarray*}
Le dernier terme est présent dans l'estimation de la norme $\|\omega^{0,m}_{\eta_{[m]} }\|_{r_1, h+1}$ seulement si $m=1$, car 
$\omega^{0,0}_{\eta_{[m]} }=\omega^{0,0} $ si $m\geq 2$. En utilisant l'hypothèse faite sur le poids $S_{h+1}$ et l'estimation (11) on a:
 \begin{eqnarray*}
\|\omega^{0,m}_{\eta_{[m]} }\|_{r_1, h+1} \leq C\cdot \sigma_m^{-s(h)}\cdot a_h^2 + 2\|\omega^{m,0}\|_{r, h}
\cdot\| T_{r,m} \,\omega^{0,m} \|_{r_1, h+1}+
\\
\\
+2\|T_{r,m}\,\omega^{0,m}\|_{r_1, h+1}\cdot\| \omega^{0,0}\|_{r, h} +c\cdot C^2\cdot \sigma ^{-2s(h)}\cdot a_h^2 \qquad\quad
\end{eqnarray*}
ce qui prouve l'estimation cherchée $\|\omega^{0,m}_{\eta_{[m]} }\|_{r_1, h+1} \leq R_m\cdot \sigma_m^{-2s(h)}\cdot a_h^2 $. Montrons maintenant l'estimation quadratique (12) pour $1\leq k<m$ (si $m\geq 2$, autrement il n'y à plus rien a prouver) en admettant qu'elle est vraie pour $k+1$. En effet en considérant les égalités suivantes on a:
\begin{eqnarray*}
\omega^{s,k}_{\eta_{[k]} }= \bar{\partial}_{_J }\eta^{s,k}+  \omega^{s+k,0}\wedge \eta^{s,k}-(-1)^{k}\eta^{s,k}\wedge \omega^{s,0}_{\eta_{[k]} } +\omega^{s,k}_{\eta_{[k+1]} }=
\\
\\
=T_{r_{m-k}, k+1}\,\bar{\partial}_{_J }\omega^{s,k}_{\eta_{[k+1]} }+\omega^{s+k,0}\wedge \eta^{s,k}-(-1)^{k}\eta^{s,k}\wedge \omega^{s,0}_{\eta_{[k]} }           
\end{eqnarray*}
En utilisant la relation (3)  pour les matrices $\omega ^{\bullet,\bullet} _{\eta_{[k+1]} } $ et en explicitant les dépendances effectives de $\eta_{[k+1]} $ dans celles ci on a:
\begin{eqnarray*}
\omega^{s,k}_{\eta_{[k]} }=(-1)^k  T_{r_{m-k},k+1}\,(\omega^{s-1,k+1}_{\eta_{[k+1]} }\wedge \omega^{s,-1})+ 
 T_{r_{m-k}, k+1}\,(\omega^{s+k+1,-1}\wedge \omega^{s,k+1}_{\eta_{[k+1]} }) -
\\
\\
-(-1)^k T_{r_{m-k},k+1}\,(\omega^{s,k}_{\eta_{[k+1]} }\wedge \omega^{s,0})-
 T_{r_{m-k}, k+1}\,(\omega^{s+k,0}\wedge \omega^{s,k}_{\eta_{[k+1]} })-\qquad
\\
\\
-\sum_{j=1}^{k-1} (-1)^{k-j} T_{r_{m-k},k+1}\, (\omega^{s+j,k-j}\wedge \omega^{s,j})
+\omega^{s+k,0}\wedge \eta^{s,k}-(-1)^{k}\eta^{s,k}\wedge \omega^{s,0}_{\eta_{[k]} } 
\end{eqnarray*}
(remarquons que $\omega^{\bullet,0}_{\eta_{[k+1]} }=\omega^{\bullet,0}$ car $k\geq 1$).
On estime donc la norme des matrices $\omega^{s,k}_{\eta_{[k]} }$ à l'aide de l'expression précédente.
\begin{eqnarray*}
\|\omega^{\bullet,k}_{\eta_{[k]} }\|_{r_{m-k+1},h+1 } \leq 2c\cdot C\cdot \sigma_m^{-s(h)}\cdot\|\omega^{\bullet,k+1}_{\eta_{[k+1]}}\|_{r_{m-k},h }+ 
\qquad\qquad\qquad\quad
\\
\\
+2c\cdot  C\cdot \sigma_m^{-s(h)}(\|\omega^{\bullet,k}\|_{r_{m-k},h } +2c\|\eta^{\bullet,k+1}\|_{r_{m-k},h } )\cdot\|\omega^{\bullet,0}\|_{r_{m-k},h }
+C\cdot \sigma_m^{-s(h)}\cdot a^2_h+
\\
\\
+\|\omega^{s+k,0}\wedge \eta^{s,k} \|_{r_{m-k+1},h+1 }+
\|\eta^{s,k}\wedge \omega^{s,0} \|_{r_{m-k+1},h+1 }+\qquad\qquad\quad
\\
\\
+2c\cdot\|\eta^{s,k}  \|_{r_{m-k+1},h+1 }\cdot \|\eta^{\bullet,k } \|_{r_{m-k+1},h+1 }\qquad\qquad\qquad\qquad\quad
\end{eqnarray*}
Le dernier terme est présent dans l'inégalité précédente seulement si $k=1$, car $\omega^{\bullet,0}_{\eta_{[k]}}=\omega^{\bullet,0}$ si $k\geq 2$. En utilisant l'hypothèse inductive, l'inégalité (11) et l'hypothèse faite sur les poids $S_{h+1}$ on aura l'estimation suivante:
\begin{eqnarray*}
\|\omega^{\bullet,k}_{\eta_{[k]} }\|_{r_{m-k+1},h+1 } \leq 
2c\cdot C\cdot R_{k+1}\cdot  \sigma_m^{-[2(m-k)+1]\cdot s(h)}\cdot a^2_h+
(2c\cdot C+C)\cdot \sigma_m^{- s(h)}\cdot a^2_h +
\\
\\
+4c^2\cdot C\cdot L_{k+1}\cdot  \sigma ^{-(m-k+2)\cdot s(h)}\cdot a^2_h+2\|\omega ^{s+k,0}\|_{r,h}\cdot\|\eta^{s,k} \|_{r_{m-k+1},h+1 } +\qquad\quad 
\\
\\
+2\|\eta^{s,k}\|_{r_{m-k+1},h+1 } \cdot \|\omega ^{s,0}\|_{r,h}
+2c\cdot L_k^2 \cdot \sigma ^{-2(m-k+1)\cdot s(h)}\cdot a^2_h \qquad\qquad\quad
\end{eqnarray*}
 ce qui prouve l'estimation (12) pour tous les indices $k=1,...,m$. Venons-en maintenant à la preuve de l'estimation (12) pour $k=0$ (et $m\geq 0$). On pose par définition
$\theta^{s,0}:= (\I_{p_s}+\eta^{s,0})^{-1}-\I_{p_s}\in M_{p_s,p_s}({\cal E}_X(B_{r_m})) $. On peut alors écrire 
$\omega ^{s,0}_{\eta_{[0]} }\equiv\omega ^{s,0}_{\eta} $ sous la forme
\begin{eqnarray*}
\omega ^{s,0}_{\eta_{[0]} }=g^{-1}_s\cdot\Big(\bar{\partial}_{_J }\eta^{s,0}+\omega ^{s,0}\wedge \eta^{s,0}+
\eta^{s-1,1}\wedge \theta ^{s-1,0}\wedge \omega ^{s,-1}+ \eta^{s-1,1}\wedge g^{-1}_{s-1}\cdot \omega ^{s,-1}\wedge \eta^{s,0}+\omega ^{s,0}_{\eta_{[1]} }\,  \Big)=
\\
\\
=g^{-1}_s\cdot\Big(T_{r_m, 1} \,\bar{\partial}_{_J }\omega ^{s,0}_{\eta_{[1]} }   +\omega ^{s,0}\wedge \eta^{s,0}+
\eta^{s-1,1}\wedge \theta ^{s-1,0}\wedge \omega ^{s,-1}+ \eta^{s-1,1}\wedge g^{-1}_{s-1}\cdot \omega ^{s,-1}\wedge \eta^{s,0} \Big)=\;
\end{eqnarray*}
En utilisant la relation (3) pour la matrice $\omega ^{s,0}_{\eta_{[1]} }$ et en remarquant que $\omega ^{\bullet,-1}_{\eta_{[1]} }=\omega ^{\bullet,-1}$ on obtient l'expression
\begin{eqnarray*}
\omega ^{s,0}_{\eta_{[0]} }
=g^{-1}_s\cdot\Big(-T_{r_m, 1}(\omega ^{s,0}_{\eta_{[1]}}\wedge\omega ^{s,0}_{\eta_{[1]}}) + T_{r_m, 1}(\omega ^{s-1,1}_{\eta_{[1]}}\wedge\omega ^{s,-1}) 
+ T_{r_m, 1}(\omega ^{s+1,1}\wedge \omega ^{s,1}_{\eta_{[1]}} ) 
\\
\\
+\omega ^{s,0}\wedge \eta^{s,0}+
\eta^{s-1,1}\wedge \theta ^{s-1,0}\wedge \omega ^{s,-1}+ \eta^{s-1,1}\wedge g^{-1}_{s-1}\cdot \omega ^{s,-1}\wedge \eta^{s,0} \Big)
\end{eqnarray*}
Le fait que $\sup_{z\in B_{r(1-\sigma )}}\|\eta^{s,0}(z)\| <\varepsilon <1$ implique l'égalité
$$
\theta ^{s,0} ={\displaystyle \sum_{j=1}^{\infty}\,(-1)^j(\eta^{s,0})^j}  
$$
(a priori la série précédente est convergente en topologie ${\cal C}^{0}$ vers l'élément $\theta ^{s,0}$ de classe ${\cal C}^{\infty}$, mais une étude élémentaire plus précise, dont on n'aura pas besoin ici, montre que l'estimation précédente $\sup_{z\in B_{r(1-\sigma )}}\|\eta^{s,0}(z)\| <\varepsilon <1$ est suffisante pour assurer la convergence de la série en topologie ${\cal C}^h$ pour tout $h$). L'inégalité
$$
\|\eta^{s,0}\|_{r(1-\sigma ),h+1}\leq L_0\cdot \sigma _m^{-(m+1)\cdot s(h)}\cdot a_h<\varepsilon <1/2
$$ 
implique la convergence de la série précédente en topologie ${\cal C}^{h+1,\mu }(\bar{B}_{r(1-\sigma )})$ et permet d'effectuer les estimations 
$$
\|\theta^{s,0}\|_{r(1-\sigma ),h+1}\leq \sum_{j=1}^{\infty}\|\eta^{s,0}\|_{r(1-\sigma ),h+1}^j\leq 2\|\eta^{s,0}\|_{r(1-\sigma ),h+1}
$$
et $\|g^{\pm 1}_{\bullet}  \|_{r(1-\sigma ),h+1}<2 $. En utilisant alors l'estimation (11) et l'estimation (12) pour $k=1$, prouvée précédemment on obtient  les estimations
\begin{eqnarray*}
\|\omega^{s,0}_{\eta_{[0]} }\|_{r(1-\sigma ),h+1} \leq 2C\cdot \sigma_m ^{-s(h)}\cdot (a_h +2L_1\cdot \sigma_m ^{-m\cdot s(h)}\cdot a_h)^2+
4C\cdot \sigma_m ^{-s(h)}\cdot c\cdot R_1\cdot \sigma_m ^{-2m\cdot s(h)}\cdot a^2_h+
\\
\\
+4\|\omega ^{s,0}\|_{r,h}\cdot\|\eta^{s,0}\|_{r(1-\sigma ),h+1}+
2c\|\eta^{s-1,1}\|_{r(1-\sigma ),h+1}\cdot \|\theta ^{s-1,0}\|_{r(1-\sigma ),h+1}+\qquad\qquad
\\
\\
+4c \|\eta^{s-1,1}\|_{r(1-\sigma ),h+1}\cdot\|\eta^{s,0}\|_{r(1-\sigma ),h+1}\leq \qquad\qquad\qquad\qquad\qquad\qquad
\end{eqnarray*}
\begin{eqnarray*}
\leq 2C\cdot \sigma_m ^{-s(h)}\cdot a^2_h +8C\cdot L_1^2\cdot \sigma_m ^{-(2m+1)\cdot s(h)}\cdot  a^2_h+4C\cdot L_1\cdot \sigma_m ^{-(m+1)\cdot s(h)}\cdot  a^2_h+\qquad\qquad
\\
\\
+4c\cdot C\cdot R_1\cdot \sigma_m ^{-(2m+1)\cdot s(h)}\cdot  a^2_h+4L_0\cdot \sigma_m ^{-(m+1)\cdot s(h)}\cdot  a^2_h+4c\cdot L_1\cdot L_0\cdot \sigma_m ^{-(2m+1)\cdot s(h)}\cdot  a^2_h\qquad
\end{eqnarray*}
ce qui prouve l'estimation (12) dans le cas $k=0$. On remarque aussi que les calculs faits pour montrer l'estimation (12) montrent aussi l'estimation 
\begin{eqnarray}
\|T_{r_{m-k}, k+1}\,\bar{\partial}_{_J }\omega^{s,k}_{\eta_{[k+1]} }\|_{r(1-\sigma ),h+1}\leq Q\cdot\sigma^{-b(m,k,h)}_m \cdot a_h^2 
\end{eqnarray}
pour tout $k=0,...,m$, (où $Q_k=Q_k(C,c)>0$ est une constante positive)
que on utilisera pour prouver l'estimation (10) sous la forme plus précise suivante:
\begin{eqnarray}
\|\omega^{\bullet,k}_{\eta} \|_{r(1-\sigma ),h+1} \leq R'_k \cdot\sigma_m ^{-[(k+2)\cdot m+1]\cdot s(h)}\cdot a_h^2
\end{eqnarray}
 où $R'_k=R'_k(C,c)>0$ est une constante positive. Considérons pourtant l'expression suivante de $ \omega^{s,k}_{\eta}$, pour $k=0,...,m$
\begin{eqnarray*}
\omega^{s,k}_{\eta}= g_{s+k}^{-1}\cdot \Big(  \bar{\partial}_{_J }\eta^{s,k}+\sum_{j=0}^{k}  \omega^{s+j,k-j}\wedge \eta^{s,j}-\sum_{j=0}^{k-1} (-1)^{k-j}\eta^{s+j,k-j}\wedge \omega^{s,j}_{\eta} +\qquad\;
\\
\\
+(-1)^k\eta^{s-1,k+1}\wedge \theta ^{s-1,0}\wedge \omega ^{s,-1}+(-1)^k\eta^{s-1,k+1}\wedge g_{s-1}^{-1}\cdot\omega ^{s,-1}\wedge \eta^{s,0}+\omega^{s,k}_{\eta_{[k+1]} }\,\Big)=
\\
\\
= g_{s+k}^{-1}\cdot \Big( T_{r_{m-k},k+1}\,  \bar{\partial}_{_J }\omega^{s,k}_{\eta_{[k+1]} } +\sum_{j=0}^{k}  \omega^{s+j,k-j}\wedge \eta^{s,j}-\sum_{j=0}^{k-1} (-1)^{k-j}\eta^{s+j,k-j}\wedge \omega^{s,j}_{\eta} +
\\
\\
+(-1)^k\eta^{s-1,k+1}\wedge \theta ^{s-1,0}\wedge \omega ^{s,-1}+(-1)^k\eta^{s-1,k+1}\wedge g_{s-1}^{-1}\cdot\omega ^{s,-1}\wedge \eta^{s,0}\Big) \quad\qquad
\end{eqnarray*}
Dans le cas $k=0$ l'estimation (14) dérive banalement de l'estimation (12) étant $\omega^{s,0}_{\eta}\equiv \omega^{s,0}_{\eta_{[0]}}$. On montre  maintenant l'estimation quadratique (14), pour tous les indices à l'aide d'une récurrence croissante sur $k=0,...,m-s$, appliquée à l'expression précédente de la matrice $\omega^{s,k}_{\eta}$ et à l'aide de l'estimation quadratique (13) . On suppose vraie l'estimation (14) pour tout $j=0,...,k-1$ et on considère l'estimation suivante en rappelant l'hypothèse faite sur le poids $S_{h+1}$ de la norme de Hölder et l'inégalité $a_h\leq 1$:
\begin{eqnarray*}
\|\omega^{s,k}_{\eta}\|_{r(1-\sigma) ,h+1} \leq 2Q_k\cdot \sigma_m ^{-2(m-k+1)\cdot s(h)}\cdot a^2_h +
4\sum_{j=0}^k  \|\omega^{s+j,k-j}\|_{r,h}\cdot \| \eta^{s,j}\|_{r(1-\sigma ),h+1}+ \quad
\\
\\
+\sum_{j=0}^{k-1}\| \eta^{s+j,k-j} \|_{r(1-\sigma ),h+1}\cdot\|  \omega^{s,j}_{\eta}\|_{r(1-\sigma ),h+1}
+8c\| \eta^{s-1,k+1}\|_{r(1-\sigma ),h+1} \cdot\| \eta^{\bullet,0}\|_{r(1-\sigma ) ,h+1}\leq 
\\
\\
\leq (2q_k+4(k+1)\cdot L)\cdot \sigma_m ^{-2(m+1)\cdot s(h)}\cdot a^2_h +\sum_{j=0}^{k-1}(L_{k-j}\cdot R'_j) \cdot \sigma_m ^{-[m+(j+2)\cdot m+1]\cdot s(h)}\cdot a^2_h+\;
\\
\\
+8c\cdot L_{k+1}\cdot L_0 \cdot \sigma _m^{-(2m+1)\cdot m\cdot s(h)}\cdot a^2_h \qquad\qquad\qquad\qquad\qquad\qquad
\end{eqnarray*}
ce qui prouve l'estimation (14) et donc l'estimation (10). \hfill $\Box$

%%%%%%%%%%%%%%%%%%%%%%%%%%%%%%%%%%%%%%%%%%%%%%%%%%%%%%
\subsubsection{Bon fonctionnement du procédé itératif.}
%%%%%%%%%%%%%%%%%%%%%%%%%%%%%%%%%%%%%%%%%%%%%%%%%%%%%%

Dans cette partie on va établir les hypothèses qui permettent d'appliquer l'estimation fondamentale (10) à une étape quelconque du procédé itératif. On commence par préciser les paramètres qui contrôleront les étapes de la convergence rapide. Pour cela on définit d'abord les paramètres $\sigma _k:=e^{-k-2}$ qui contrôleront les restrictions des rayons des boules, lesquels sont définis de façon récursive par la formule $r_{k+1}:=r_k(1-\sigma _k)$ pour tout entier $k \geq 0$. Bien évidement le rayon limite 
$$
r_{\infty}:=\lim_{k\rightarrow +\infty}r_k =r_0{\displaystyle \prod_{k=0}^{\infty}(1-\sigma _k)}
$$
est non nul, car  $-\sum_{k=0}^{\infty}\log(1-\sigma _k)\leq Cst \sum_{k=0}^{\infty}\sigma _k <\infty$. 
Ensuite on désigne par $r(k,l):=r_k(1-l \cdot \sigma _{m,k}),\;l=1,...,m+1$ où $\sigma _{m,k}:=\sigma _k/(m+1)$. Pour le choix du rayon initial on considère les suites numériques
\begin{eqnarray*} 
\alpha _k(r):= a_0(r)^{2^k}\cdot \prod_{j=0}^{k-1} H^{2^{k-1-j}}\cdot\sigma _j^{-\nu (m,j)\cdot2^{k-1-j}}
\\
\\   
\beta  _k(r):= b_0(r)^{2^k}\cdot \prod_{j=0}^{k-1} P^{2^{k-1-j}}\cdot e^{-\gamma (m,j)\cdot 2^{k-1-j}} 
\end{eqnarray*}
où $b_0(r):=H\cdot \sigma _{m,0}^{-(m+1)\cdot s(0) }\cdot a_0(r)$, $H>0,\;P>0$ sont deux constantes (dépendants seulement de $m$) et $\gamma (m,j) $ est une fonction affine strictement croissante en $j$.  On rappel que par définition de la norme de Hölder on a que la quantité $a_0=a_0(r)$ tend vers zéro lorsque $r>0$ tend vers zéro. Avant de faire le choix du rayon initial on a besoin du lemme élémentaire suivant.
\addtocounter{subsubsection}{-1}  
\begin{lem}. Il existe $\rho \in (0,1)$ tel que;
\\
$(A)$, pour tout $r\in (0,\rho ]$  les séries numériques $\sum_{k\geq 0}\alpha _k(r) $ et $\sum_{k\geq 0}\beta_k(r)$  sont convergentes.
\\
$(B)$
\begin{eqnarray*}
\lim_{r\rightarrow 0^+ } \sum_{k=0}^{\infty}\,\alpha _k(r)=0, \quad\lim_{r\rightarrow 0^+ } \sum_{k=0}^{\infty}\,\beta _k(r)=0 
\end{eqnarray*}
\\
$(C)$, pour tout $r\in (0,\rho ]$ et $K\geq 0$ on a les inégalités $\alpha _k(r)\leq 1$ et $\beta _k(r)<\varepsilon$
\end{lem}
$Preuve$. Par le critère du rapport il suffit de vérifier que les suites $\ln(\alpha _{k+1}/ \alpha _k)$ et $\ln(\beta  _{k+1}/\beta_k)$ tendent vers $-\infty$ lorsque $k$ tend vers $+\infty$. En explicitant les logarithmes on a:
\begin{eqnarray*} 
\ln(\alpha _{k+1}/ \alpha _k )=2^k\cdot\Big( \ln a_0(\rho)+2^{-1}(\ln H)\sum_{j=0}^{k-1}2^{-j}-2^{-1}\sum_{j=0}^{k-1}\nu (m,j)2^{-j}\ln \sigma _{m,j}  \Big)-\qquad
\\
\\
-\nu (m,k)\ln\sigma _{m,k} +\ln H\qquad
\\
\\
\ln(\beta  _{k+1}/\beta  _k) =2^k\cdot\Big( \ln b_0(\rho )+2^{-1}(\ln P)\sum_{j=0}^{k-1}2^{-j} + 2^{-1}\sum_{j=0}^{k-1}\gamma (m,j)2^{-j}\Big)+\gamma (m,k)+\ln P
\end{eqnarray*}
Il suffit donc choisir $\rho >0$ suffisamment petit pour assurer les inégalités
\begin{eqnarray*} 
 \ln a_0(\rho )+2^{-1}(\ln H)\sum_{j=0}^{\infty}2^{-j}+2^{-1}\sum_{j=0}^{\infty}\nu (m,j)[j+2+\ln (m+1)]2^{-j}<0
\\
\\
\ln b_0(\rho )+2^{-1}(\ln P)\sum_{j=0}^{\infty}2^{-j} +  2^{-1}\sum_{j=0}^{\infty}\gamma (m,j)2^{-j}<0\qquad\qquad
\end{eqnarray*}
(se rappeler la définition de $\sigma _j$). On aura alors la convergence voulue pour les suites $\ln(\alpha _{k+1}/ \alpha _k)$ et $\ln(\beta  _{k+1}/\beta_k)$. Le fait que pour tout entier $k\geq 0$ les quantités $\alpha_k(r)$ et $\beta_k(r)$ tendent monotonement vers zéro lorsque $r$ tend vers zéro, implique par le théorème de la convergence dominée, les conclusions $(B)$ et $(C)$ du lemme.\hfill $\Box$
\\
D'après le lemme précédent  on peut  choisir le rayon initial $r_0\in (0,\rho ]$ de telle sorte que l'inégalité $\sum_{k=0}^{\infty}\,\beta _k(r_0)<1/2 + 1/(4\ln 2)$ soit satisfaite. Dans la suite on notera $\alpha _k:=\alpha _k(r_0)$ et $\beta _k:=\beta _k(r_0)$.
On définit maintenant le paramètre $\eta_{k+1}$, qui contrôle la recalibration des éléments $\omega _k$, $k\geq 0$ au $k$-ème pas du procédé itératif (on désigne avec $\omega _0=\omega$ le choix initial de $\omega$ associée au système $(S)=(S_{\omega})$), de la façon suivante; on définit récursivement en ordre décroissant sur $t=m,...,0$ les matrices 
$$
\eta^{s,t}_{k+1} :=-T_{r(k,m-t),t} \,\Big(\omega ^{s,t}_k+\omega ^{s+t+1,-1}_k\wedge\eta^{s,t+1}_{k+1} +(-1)^t\eta^{s-1,t+1}_{k+1}
\wedge\omega ^{s,-1}_k\Big)\in M_{p_{s+t},p_s}({\cal E}_{X}^{0,t}(B_{r(k,m-t)}))     
$$
et on pose $\eta_{k+1} :=(\eta^{s,t}_{k+1})_{s,t}\in {\cal P}(\bar{B}_{r(k,m)}) $. On va justifier ensuite l'estimation $\sup_{z\in B_r} \|\eta^{s,0} _{k+1} (z)\| <\varepsilon $ qui assure l'invertibilité de la matrice $g_{s,k+1}$. On définit alors les matrices
$\omega^{s,t}_{k+1}:=\omega^{s,t}_{k,\;\eta_{k+1}}\in M_{p_{s+k},p_s}({\cal E}_{X}^{0,t+1}(\bar{B}_{r_{k+1} }))$ 
pour tout $t=-1,...,m$. Les constantes $R_k$, $L_k$ qui apparaissent dans la définition des poids $S=(S_k)_{k\geq 0}$  de la norme de Hölder sont définies par les formules:
\begin{eqnarray*}
R_{k+1} :=\max\{\|\omega ^{s,t}_{k,\,I} \|_{r_k,\,\mu }\,/\|\partial^{k+1} \omega^{s,t}_{k,\,I}  \|_{r_k ,\,\mu }  
 \,|\,0\leq s \leq m,\,0\leq t \leq m-s,\,|I|=t+1\}       
\\
\\
L_{k+1} :=\max\{ 2^{-k-1} \|g_s(k)^{\pm 1}   \|_{r_k ,\,\mu }\,/\|\partial^{k+1} g_s(k)^{\pm 1}  \|_{r_k ,\,\mu }  \,|\,0\leq s \leq m\}  \qquad\qquad
\end{eqnarray*}
pour tout entier $k\geq 0$. Ici on suppose que 
$\max\{\|\omega ^{s,t}_{k,\,I} \|_{r_k,\,\mu }   
 \,|\,0\leq s \leq m,\,0\leq t \leq m-s,\,|I|=t+1\}>0$, autrement il n'y a rien à prouver.
 Avec ce choix des poids $S_k>0$ on aura que les inégalités
\begin{eqnarray}
S_{k+1} \|\partial^{k+1} \omega^{\bullet,t}_{k,\,I}  \|_{r_k ,\,\mu }\leq\| \omega^{\bullet,t}_{k,\,I}  \|_{r_k ,\,\mu } \qquad\quad
\\\nonumber
\\
S_{k+1} \|\partial^{k+1} g_{\bullet}(k)^{\pm 1}  \|_{r_k ,\,\mu }\leq 2^{-k-1}\|g_{\bullet}(k)^{\pm 1}   \|_{r_k ,\,\mu }  
\end{eqnarray}
seront satisfaites pour tout  entier $k\geq 0$ et $t=0,...,m$. On pose par définition
$$
a_{k} :=\max\{\|\omega ^{s,t}_k\|_{r_k,\,k}\,|\,0\leq s \leq m\,,0\leq t \leq m-s\}
$$
$b_k:=H\cdot\sigma_{m,\,k}  ^{-(m+1)\cdot s(k) }\cdot a_k $ et $c\equiv c(\omega _0)$. Avec les notations introduites précédemment on a la proposition suivante.
\begin{prop} 
Pour tout entier $k\geq 0$ on a les estimations suivantes;
\begin{eqnarray}
a_{k+1} \leq H\cdot\sigma_{m,\,k} ^{-\nu (m,k)}\cdot a^2_{k}\leq 1\;
\\\nonumber
\\
\|\eta^{s,t}_{k+1} \|_{r_{k+1}  ,\,k+1 } \leq  b_k   <\varepsilon <1/2
\\\nonumber
\\
\|\omega ^{s,-1} _{k+1}\|_{r_{k+1},k+1 } \leq 4c    \qquad
\end{eqnarray}
où $H:=\max\{R(C,4c), L(C,4c)\}>0$ et les inégalités $a_k\leq \alpha _k,\;b_k\leq \beta _k$.
\end{prop} 
$Preuve$. On montre les trois estimations précédentes à l'aide d'une récurrence sur $k\geq 0$. Pour $k=0$ on a d'après le lemme (3.4.2), la validité des inégalités
$\alpha _0(r)\leq 1$ et $\beta _0(r)<\varepsilon <1/2$. On est donc en position d'appliquer la proposition (3.4), laquelle assure les inégalités (17) et (18) pour $k=0$. On obtient alors l'estimation $\|\eta^{\bullet,0}_1 \|_{r_1,1}<1/2$ laquelle, pour les calculs faits dans la preuve de la proposition (3.4), assure les inégalités 
$\|g_{\bullet,1}^{\pm 1} \|_{r_1,1}<2$ qui assurent donc l'estimation (20) pour $k=0$, car 
$$
\|\omega ^{s,-1} _1\|_{r_1,1} \leq \|g_{s-1,1}^{-1} \|_{r_1,1}\cdot \|\omega ^{s,-1}\|_{r_1,1}\cdot \|g_{s,1}\|_{r_1,1}
$$ 
Supposons maintenant par hypothèse récursive que les estimations (17), (18) et (19), sont vraies pour tout $l=0,...,k-1$. L'inégalité  (17) implique alors l'inégalité 
$$
b_{l+1}\leq P\cdot e^{\gamma (m,l)}\cdot b_l^2
$$
pour les indices en question. On obtient donc les inégalités $a_l\leq \alpha _l\leq 1,\;b_l\leq \beta _l<\varepsilon $. Le fait que par hypothèse inductive on dispose de l'estimation $\|\omega ^{\bullet,-1}_k\|_{r_k,\,k}<4c$ 
permet alors d'appliquer la proposition (3.4) laquelle assure la validité des estimations (17) et (18) pour $l=k$. En particulier l'estimation (18) assure la validité de l'estimation
$$
\|\eta^{\bullet,0}_{k+1} \|_{r_{k+1} ,\,k+1 } \leq \beta _k<\varepsilon <1/2
$$
On passe  maintenant à l'estimation des  normes $\|\omega^{s,-1}_{k+1} \|_{r_{k+1} ,\,k+1}$. On a par définition des matrices $\omega ^{s,-1}_{k+1}$ l'estimation  
$$
\|\omega ^{s,-1}_{k+1} \|_{r_{k+1} ,\,k+1} \leq \|g_{s-1}(k+1)^{-1} \|_{r_{k+1} ,\,k+1}\cdot \|\omega ^{s,-1}\|_{r_{k+1} ,\,k+1}\cdot \|g_s(k+1)\|_{r_{k+1} ,\,k+1}
$$
La condition (16) sur les poids implique que pour tout $k\geq 0$ on dispose de l'inégalité 
$$
\|g_s(k+1)^{\pm 1}\|_{r_{k+1} ,\,k+1}\leq (1+2^{-k-1})\cdot \|g_{s,k+1}^{\pm 1}\|_{r_{k+1} ,\,k+1}\cdot\|g_s(k)^{\pm 1}\|_{r_k ,\,k}
$$
L' hypothèse inductive nous permet alors d'effectuer les estimations suivantes pour $k\geq 1$
\begin{eqnarray*} 
\|g_{s-1}(k+1) \|_{r_{k+1} ,\,k+1}\leq \prod_{j=2}^{k+1}(1+2^{-j})\cdot \prod_{j=1}^{k+1}(1+\|\eta^{s-1,0}_j\|_{r_j,\,j}) \leq
\\
\\ 
\leq \sqrt{e}\cdot  \exp\Big(\sum_{j=1}^{k+1}\,\|\eta^{s-1,0}\|_{r_j,\,j}  \Big)\leq  \sqrt{e}\cdot  \exp\Big(\sum_{j=0}^{\infty}\,\beta _j \Big)<2\quad
\end{eqnarray*}
et
\begin{eqnarray*}  
\|g_s(k+1)^{-1}  \|_{r_{k+1} ,\,k+1}\leq \prod_{j=2}^{k+1}(1+2^{-j})\cdot  \prod_{j=1}^{k+1}\|g_{s,j}^{-1} \|_{r_j,\,j} \leq 
\sqrt{e}\cdot\prod_{j=1}^{k+1}(1- \|\eta^{s,0}_j\|_{r_j,\,j} )^{-1}  \leq 
\\
\\
\leq \sqrt{e}\cdot\exp\Big(2(\ln 2)\sum_{j=1}^{k+1}\|\eta^{s,0}_j\|_{r_j,\,j}\Big) \leq  \sqrt{e}\cdot \exp\Big(2(\ln 2)\sum_{j=0}^{\infty}\,\beta _j\Big)<2\qquad\qquad
\end{eqnarray*}
(rappeler le choix initial du rayon $r_0$). On obtient donc l'estimation voulue $\|\omega ^{s,-1}_{k+1} \|_{r_{k+1} ,\,k+1} \leq 4c$, ce qui conclu la preuve des trois estimations (17), (18) et (19) pour $j=k$ et donc pour tout entier positif $k$.\hfill $\Box$

%%%%%%%%%%%%%%%%%%%%%%%%%%%%%%%%%%%%%%%%%%%%%%%%%%%%%%%%%%%%%%%%%%%%%%%%%%%%%%%
\subsubsection{Convergence vers une solution du problème différentiel (S).}
%%%%%%%%%%%%%%%%%%%%%%%%%%%%%%%%%%%%%%%%%%%%%%%%%%%%%%%%%%%%%%%%%%%%%%%%%%%%%%% 
Avec les notations introduites précédemment on a la proposition suivante.
\begin{prop}. Les limites  
$$
\eta^{s,t}:=\lim_{k\rightarrow \infty} \eta (k)^{s,t}
$$
$s=0,...,m,\;t=0,...,m-s$ existent en topologie ${\cal C}^{h,\mu}(B_{r_{\infty}})$ pour tout $h\geq 0$ et ils constituent les composantes du paramètre de recalibration $\eta=(\eta^{s,t})_{s,t}\in {\cal P}(B_{r_{\infty}})$, solution du problème différentiel $(S)$.
\end{prop} 
$Preuve$. Nous commençons par prouver l'existence des limites
$$
g_s:={\displaystyle \prod_{j\geq 1 }^{\longrightarrow}}g_{s,j}=\lim _{k\rightarrow \infty}g_s(k)=\I_{p_s}+\lim_{k\rightarrow \infty} \eta (k)^{s,0} 
$$
 en topologie ${\cal C}^{h,\mu}(B_{r_{\infty}})$ pour  $h\geq 0$ quelconque et le fait que les matrices $g_s$ sont inversibles. On aura alors $\eta^{s,0}=g_s-\I_{p_s}$. On déduit immédiatement de la proposition (3.5), les estimations 
$\|\eta^{\bullet,\bullet}_k\|_{r_{\infty},h }<\beta _k<1/2$
et 
$\|g_{\bullet} (k)^{\pm 1}  \|_{r_{\infty},h }<2$ pour tout $k\geq h$. Le fait que 
$g_s(k+1)-g_s(k)=\eta^{s,0}_{k+1}\cdot g_s(k)$ implique les estimations suivantes:
$$
\|g_s(k+1)-g_s(k)\|_{r_{\infty} ,\,h}\leq \|\eta^{s,0}_{k+1}\|_{r_{\infty} ,\,h}\cdot \|g_s(k)\|_{r_{\infty} ,\,h}\leq 2\|\eta^{s,0} _{k+1}\|_{r_{\infty} ,\,h}\leq 2\beta _k
$$
pour tout $k\geq h$  
On a alors
\begin{eqnarray*}
\sum_{k=0}^{\infty}\|g_s(k+1)-g_s(k)\|_{r_{\infty} ,\,h}\leq\sum_{k=0}^h\|g_s(k+1)-g_s(k)\|_{r_{\infty} ,\,h}+2\sum_{k=h+1}^{\infty}\,\beta _k <\infty 
\end{eqnarray*}
et donc l'existence des matrices $g_s\in{\cal C}^{h,\mu}(B_{r_{\infty}} ,M_{p_s,p_s}(\C))$ pour tout $h \geq 0$  telles que 
\begin{eqnarray*}
g_s=\lim_{k\rightarrow \infty}g_s(k)=\I_{p_s}+  \sum_{k=0}^{\infty}(g_s(k+1)-g_s(k))
\end{eqnarray*}
en topologie ${\cal C}^{h,\mu}(B_{r_{\infty}})$. D'autre part l'égalité
\begin{eqnarray*}
g_s(k+1)^{-1} -g_s(k)^{-1} =g_s(k)^{-1}\sum_{j=1}^{\infty} \,(-1)^j (\eta^{s,0}_{k+1})^j  
\end{eqnarray*}
permet d'effectuer les estimations suivantes
\begin{eqnarray*}
\|g_s(k+1)^{-1}-g_s(k)^{-1}\|_{r_{\infty} ,\,h}\leq \|g_s(k)^{-1}\|_{r_{\infty} ,\,h}\cdot \sum_{j=1}^{\infty}\|\eta^{s,0}_{k+1} \|_{r_{\infty} ,\,h}^j
\leq 4\|\eta^{s,0} _{k+1} \|_{r_{\infty} ,\,h}\leq 4\beta _k
\end{eqnarray*}
pour tout $k\geq h$, lesquelles  assurent la convergence de la série
\begin{eqnarray*}
\sum_{k=0}^{\infty}\|g_s(k+1)^{-1} -g_s(k)^{-1} \|_{r_{\infty} ,\,h}\leq
\sum_{k=0}^h\|g_s(k+1)^{-1} -g_s(k)^{-1} \|_{r_{\infty} ,\,h} + 4\sum_{k=h+1}^{\infty} \beta_k<\infty 
\end{eqnarray*}
ce qui prouve l'existence des matrices $\rho_s\in{\cal C}^{h,\mu}(B_{r_{\infty} } ,M_{p_s,p_s}(\C))$ telles que $\rho _s=\lim_{k\rightarrow \infty}g_s(k)^{-1}$ en topologie ${\cal C}^{h,\mu}(B_{r_{\infty}})$ pour tout $h\geq 0$. On a alors l'égalité  $\I_{p_s}=\lim_{k\rightarrow \infty}g_s(k)\cdot g_s(k)^{-1}=g_s\cdot \rho _s$, qui montre que 
$g_s\in GL(p_s, {\cal E}(B_{r_{\infty}}))$ et $\rho _s=g^{-1}_s$. Montrons maintenant l'existence des limites $\eta^{s,t}$ pour $t\geq 1$. En rappelant l'expression des composantes $\eta(k)^{s,t},\;t\geq 1$, du paramètre $\eta(k)$ introduite dans la sub-section (3.2) et en tenant compte de l'existence de la limite 
$g\in \Gamma (U)$ prouvé précédemment, on déduit que il suffit de prouver l'existence des limites
\begin{eqnarray*}
\lim_{k\rightarrow \infty } \sum_{J\in J_k(\rho (\tau) )}\;{\displaystyle \bigwedge_{1\leq r\leq \rho (\tau) }^{\longrightarrow}} 
\;g_{s+\sigma' (\tau,r)}  (j_r-1)\cdot \eta^{s+\sigma (\tau,r),\,\tau_{\rho (\tau)+1-r } }_{j_r} \cdot g_{s+\sigma (\tau,r)} (j_r)^{-1}
\end{eqnarray*}
 $\tau\in \Delta _t,\;t\geq 1$ en topologie ${\cal C}^{h,\mu}(B_{r_{\infty} })$, avec $h\geq 0$ quelconque, pour obtenir l'existence de la limite 
$\eta^{s,t}\in M_{p_{s+t},p_s}({\cal E}_{_X}^{0,t}(B_{r_{\infty} }))$. Il suffit donc de prouver que pour tout $h\geq \rho (\tau)$ la limite 
\begin{eqnarray*}
\lim_{k\rightarrow \infty } \sum_{J\in J_k(\rho (\tau) )}\;\prod_{r=1}^{\rho (\tau)}   
\;\|g_{\bullet}  (j_r-1)\cdot \eta^{\bullet,\bullet}_{j_r} \cdot g_{\bullet} (j_r)^{-1}\|_{r_{\infty},h} =\qquad
\\
\\
=\lim_{k\rightarrow \infty }{\displaystyle\sum_{\scriptstyle l+p=\rho (\tau) 
\atop
\scriptstyle  l,p\geq 0} }\;\sum_{I\in J_h(l)}\;\prod_{r=1}^l 
\;\|g_{\bullet} (i_r-1)\cdot \eta^{\bullet,\bullet}_{i_r} \cdot g_{\bullet} (i_r)^{-1}\|_{r_{\infty},h}\times
\\
\\
\times\sum_{I\in J_{h,k} (p)}\;\prod_{r=1}^p 
\;\|g_{\bullet} (j_r-1)\cdot \eta^{\bullet,\bullet}_{j_r} \cdot g_{\bullet} (j_r)^{-1}\|_{r_{\infty},h}\qquad\quad
\end{eqnarray*}
est finie. Ici on pose par définition $J_{h,k} (p):=\{J\in \{h+1,...,k\}^p\;|\;j_1<...<j_p \},\;J_h(0)=J_{h,k}(0):=\emptyset$, (rappeler la convention faite sur les symboles de somme et de produit). On considère pourtant les estimations suivantes pour $1\leq p\leq \rho (\tau)$
\begin{eqnarray*}
\lim_{k\rightarrow \infty }\sum_{I\in J_{h,k} (p)}\;\prod_{r=1}^p 
\;\|g_{\bullet} (j_r-1)\|_{r_{\infty},h}\cdot \|\eta^{\bullet,\bullet}_{j_r}\|_{r_{\infty},h} \cdot \|g_{\bullet} (j_r)^{-1}\|_{r_{\infty},h}<\qquad\qquad\quad
\\
\\
<4^p\lim_{k\rightarrow \infty }\sum_{J\in J_{h,k} (p)}\;\prod_{r=1}^p \beta_{j_r-1}
<4^p\lim_{k\rightarrow \infty }\sum_{J\in J_{h,k} (p)}\beta_{j_p-1}
=4^p\lim_{k\rightarrow \infty }\sum_{j=h+p}^k \,{j-h-1\choose p-1 }\cdot \beta _{j-1}  <\infty    
\end{eqnarray*}
La dernière limite est finie par le critère du rapport, rappeler en fait que $\lim_{k\rightarrow \infty }\beta _{k+1}/\beta _k=0$, d'après la preuve du lemme 
(3.4.2). On a donc prouvé l'existence du paramètre limite $\eta \in {\cal P}(B_{r_{\infty} })$. Prouvons maintenant qu'il constitue une solution 
(de classe ${\cal C}^{\infty}$) pour le problème différentiel (S). En effet l'existence de la limite $\eta$ en topologie 
${\cal C}^{h,\mu}(B_{r_{\infty} })$, pour $h\geq 1$ implique les égalités
\begin{eqnarray*}
\omega ^{s,t}_{\eta}=\lim_{k\rightarrow \infty}\omega ^{s,t}_{\eta(k)} =\lim_{k\rightarrow \infty}\omega ^{s,t}_k
\end{eqnarray*} 
en topologie ${\cal C}^{h-1,\mu}(B_{r_{\infty} })$. En rappelant l'inégalité $a_k\leq \alpha _k$ obtenue dans la preuve de la proposition (3.5), on obtient
\begin{eqnarray*}
\|\omega ^{s,t}_{\eta} \|_{r,0}=\lim_{k\rightarrow \infty}\|\omega ^{s,t}_k\|_{r,0} \leq\lim_{k\rightarrow \infty}a_k=0  
\end{eqnarray*} 
ce qui prouve que $\eta \in {\cal P}(B_{r_{\infty} })$ est une solution du  système différentiel
$$
\left  \{
\begin{array}{lr}
\omega^{s,t}_{\eta}=0
\\
s=0,...,m
\\
t=0,...,m-s 
\end{array}                                                                           
\right.
$$
qui n'est rien d'autre que le système différentiel (S).  \hfill $\Box$

%%%%%%%%%%%%%%%%%%%%%%%%%%%%%%%%%%%%%%%%%%%%%%%%%%%%%%%%%%%%%
\subsection{Cinquième étape: fin de la preuve du théorème 1.}
%%%%%%%%%%%%%%%%%%%%%%%%%%%%%%%%%%%%%%%%%%%%%%%%%%%%%%%%%%%%%

L'étape précédente montre que on peut se ramener à considérer le diagramme commutatif suivant. 
\begin{diagram}[height=1cm,width=1cm]
&         &                              &                 &0&                                         &0                                                       \\
&         &                              &              &\uTo&                                      &\uTo                                                       \\
0&  \rTo  &(Ker\bar{\partial})_{_{|V} } &  \rTo   &{\cal G}_{|_V}& \rTo^{\bar{\partial}} 
 & {\cal G}_{|_V}\otimes_{_{{\cal E}_{_V}}}{\cal E}^{0,1}_{_V} 
\\
&          &\uTo^{\tilde{\psi}_{|..}}&                     &\uTo_{\tilde{\psi}}&                      &\uTo_{\tilde{\psi}\otimes\I_{(0,1)}}                   \\
0&\rTo  & {\cal O} _{_V}^{\oplus p_0} &  \rTo   &{\cal E}^{\oplus p_0}_{_V}&    \rTo^{\bar{\partial}_{_J }}
 & ({\cal E}^{0,1}_{_V} )^{\oplus p_0}&  
\\
&             &\uTo^{\tilde{\varphi}_{|..}} &                             &\uTo_{\tilde{\varphi}}&               &\uTo_{\tilde{\varphi}\otimes \I_{(0,1)} }
\\
0& \rTo  & {\cal O}
          _{_V}^{\oplus p_1} &  \rTo   &{\cal E}^{\oplus p_1}_{_V}&    \rTo^{\bar{\partial}_{_J }}
 & ({\cal E}^{0,1}_{_V} )^{\oplus p_1} 
\end{diagram} 
\\
avec $V:=B_{r_{\infty}}, \;\tilde{\psi}:=\psi_g, \;\tilde{\varphi}_{\bullet} =\varphi_{\bullet,\,g} $. Ce diagramme est exact, sauf pour l'instant au niveau des premières flèches verticales à gauche. Pour conclure il nous reste donc à montrer l'exactitude de la suite
$$
{\cal O}_{_V}^{\oplus p_1}\stackrel{\tilde{\varphi}_{|..}}{\longrightarrow}  
{\cal O}_{_V}^{\oplus p_0}\stackrel{\tilde{\psi}_{|..}}{\longrightarrow}
(Ker\bar{\partial})_{_{|V}}\rightarrow 0
$$
On identifie $\tilde{\varphi}=(\tilde{\varphi}_1,...,\tilde{\varphi}_{p_1}  )$,  $\tilde{\psi}=(\tilde{\psi}_1,...,\tilde{\psi}_{p_0} )$, on indique avec  
$\tilde{\varphi}_l^k \in{\cal O}_X (V)$ la $k$-ème composante de $\tilde{\varphi}_l$, et on pose 
$$
a_{k,x} :={\cal O}_x\tilde{\varphi}^k_{1,x}+...+{\cal O}_x\tilde{\varphi}^k_{p_1,x}\lhd {\cal O}_x
$$
Le fait que $(a_{k,x} \cdot {\cal E}_x)\cap {\cal O}_x=a_{k,x} $ (par définition de fidélité plate de l'anneau ${\cal E}_x$ sur l'anneau ${\cal O}_x$. On peut aussi  déduire l'égalité précédente en utilisant un résultat beaucoup plus simple, i.e la
 fidélité plate de l'anneau des séries formelles $\quad{\cal E}_x /m^{\infty}({\cal E}_x)=\hat{{\cal O}}_x$ en $x$, sur l'anneau ${\cal O}_x$ (voir\cite{mal})) implique la surjectivité du morphisme
$$
\tilde{\varphi}_{|..}:{\cal O}_{_V}^{\oplus p_1}\longrightarrow {\cal R}^{\cal O}(\tilde{\psi})
$$
où ${\cal R}^{\cal O}(\tilde{\psi})={\cal R}^{\cal E}(\tilde{\psi})\cap{\cal O}_{_V}^{\oplus p_0}$ désigne le faisceau des relations holomorphes de $\tilde{\psi}$. On pose par définition
$${\cal F}:=Im(\tilde{\psi}_{|..}: {\cal O}_{_V}^{\oplus p_1}\longrightarrow (Ker\bar{\partial})_{_{|V}})$$
L'exactitude de la suite
$
{\cal O}_{_V}^{\oplus p_1}\stackrel{\tilde{\varphi}_{|..}}{\longrightarrow}  
{\cal O}_{_V}^{\oplus p_0}\stackrel{\tilde{\psi}_{|..}}{\longrightarrow}
{\cal F} \rightarrow 0
$
et la platitude de l'anneau ${\cal E}_x$ sur l'anneau ${\cal O}_x$ impliquent l'existence du diagramme commutatif exact
\begin{diagram}[height=1cm,width=1cm]
{\cal E}^{\oplus p_1}_{_V}& \rTo^{\tilde{\varphi}_1}&{\cal E}^{\oplus p_0}_{_V}&\rTo^{{\tilde{\psi}}_{|..} \otimes \I}&{\cal F}\otimes_{_{{\cal O}_{_V}}}{\cal E}_{_V}&\rTo 0                                                                             \\
\dTo^{\I} &   &\dTo^{\I} &   &\dTo^{\alpha }&                                                                              \\
{\cal E}^{\oplus p_1}_{_V}& \rTo^{\tilde{\varphi}_1}&{\cal E}^{\oplus p_0}_{_V}&\rTo^{\tilde{\psi}}&{\cal G}_{|_V} &\rTo 0
\end{diagram} 
avec $\alpha :\sum_k\tilde{\psi}_{k,y}\otimes_{_{{\cal O}_y}} f_{k,y}\mapsto\sum_k\tilde{\psi}_{k,y} \cdot f_{k,y}\  $. On a alors 
${\cal G}_{|_V}\cong{\cal F}\otimes_{_{{\cal O}_{_V}}}{\cal E}_{_V}.$ Les égalités
\begin{eqnarray*}
\bar{\partial}_{_{\cal F}}\Big(\sum_{k=1}^{p_0}  \tilde{\psi}_{k,y}\otimes_{_{{\cal O}_y}} f_{k,y}\Big) &=&
 \sum_{k=1}^{p_0} \tilde{\psi}_{k,y}\otimes_{_{{\cal O}_y}}\bar{\partial}_{_J } f_{k,y}\cong
\sum_{k=1}^{p_0}\tilde{\psi}_{k,y}\otimes_{_{{\cal E}_y}}\bar{\partial}_{_J } f_{k,y} 
\\
\\
\bar{\partial}\Big(\sum_{k=1}^{p_0} \tilde{\psi}_{k,y}\cdot f_{k,y}\Big) & =&\sum_{k=1}^{p_0} \tilde{\psi}_{k,y}\otimes_{_{{\cal E}_y}}\bar{\partial}_{_J } f_{k,y}  
\end{eqnarray*} 
impliquent la commutativité du diagramme suivant
\begin{diagram}[height=1cm,width=1cm]
{\cal G}_{|_V}&\rTo^{\bar{\partial}} &{\cal G}_{|_V}\otimes_{_{{\cal E}_{_V}}}{\cal E}^{0,1} _{_V}&             &
\\
\uTo^{\alpha}_{\wr}&                                      &\uTo_{\wr}&          \luTo^{\alpha\otimes\I_{(0,1)}}                          \\    
{\cal F}\otimes_{_{{\cal O}_{_V}}}{\cal E}_{_V}&\rTo^{\bar{\partial}_{_{\cal F}}}&{\cal F}\otimes_{_{{\cal O}_{_V}}}{\cal E}^{0,1} _{_V}&\simeq\ 
&{\cal F}^{\infty} \otimes_{_{{\cal E}_{_V}}}{\cal E}^{0,1} _{_V}
\end{diagram} 
qui montre que 
${\cal F}=(Ker\bar{\partial})_{_{|V}},$ ce qui démontre le théorème 1.

%%%%%%%%%%%%%%%%%%%%%%%%%%%%%%%%%%%%%%%%%%%%%%%%%%%%%%%%%%%%%%%%%%%%%%%%%%%%%%%%%%%%%%%%%%%%%%%%%%%%%%%%%%%%%%%%%%%%%%%%%%%%%%%%%%%%
\section{Un résultat d'intégrabilité des connexions sur les faisceaux de ${\cal E}$-modules au dessus d'une variété différentiable.}
%%%%%%%%%%%%%%%%%%%%%%%%%%%%%%%%%%%%%%%%%%%%%%%%%%%%%%%%%%%%%%%%%%%%%%%%%%%%%%%%%%%%%%%%%%%%%%%%%%%%%%%%%%%%%%%%%%%%%%%%%%%%%%%%%%%% 

Le présent travail s'est situé principalement dans le cadre complexe car c'était notre principal intérêt. Toutefois, il est possible de  déduire aussi un résultat d'intégrabilité dans le cas des variétés ${\cal C}^{\infty}$. Considérons en effet $(X,{\cal E}_X)$ une variété ${\cal C}^{\infty}$ (${\cal E}_X\equiv{\cal E}_X(\R)$ représente ici le faisceau des fonctions ${\cal C}^{\infty}$ à valeurs réelles) et $D:{\cal G}\longrightarrow  {\cal G} \otimes_{_{{\cal E}_X }}{\cal E}(T^*_X)$ une connexion sur le faisceau de ${\cal E}_X(\K)$-modules ${\cal G}$ où $\K=\R,\;\C$. Si le faisceau des sections parallèles $KerD$ engendre ${\cal G}$ sur ${\cal E}_X(\K)$ alors évidement $D^2=0$. Le théorème suivant donne une réciproque de ce fait dans un cas particulier.
\begin{theoreme}
Soit $(X,{\cal E}_X)$ une variété différentiable et soit ${\cal G}$ un faisceau de  ${\cal E}_X(\K)$-modules,  muni d'une connexion  
$D:{\cal G}\longrightarrow  {\cal G} \otimes_{_{{\cal E}_{_X}  }}{\cal E}(T^*_X)$ telle que $D^{2}=0$. Si de plus le faisceau ${\cal G}$ admet localement une ${\cal E}(\K)$-résolution de longueur finie, alors le faisceau des sections parallèles $KerD$ engendre sur ${\cal E}_X(\K)$ le faisceau   $\cal G$ qui est le faisceau  des sections d'un système local de coefficients et le complexe $({\cal G} \otimes_{_{{\cal E}_X }}{\cal E}(\Lambda ^{\bullet} T^*_X)\; ; D)$ est une résolution acyclique du faisceau des sections parallèles. Il existe alors l'isomorphisme fonctoriel  de De Rham-Weil
$H^k(X,KerD)\cong H^k(\Gamma (X,{\cal G} \otimes_{_{ {\cal E}_X }}{\cal E}(\Lambda ^{\bullet} T^*_X))\; ; D)$.
\end{theoreme}
$Preuve$. On commence par substituire formellement dans les étapes de la preuve du théorème 1 la connexion $\bar{\partial}$ avec $D$, les opérateurs 
$\bar{\partial}_{_J }$  avec $d$ et l'opérateur de Leray-Koppelman avec l'opérateur d'homotopie de Poincaré 
$$
P_q: {\cal C}^{h} _{q+1}(B_r, M_{k,l}(\K))\longrightarrow  {\cal C}^h_q(B_r, M_{k,l}(\K))
$$
pour $q\geq 0$. Il est élémentaire de vérifier que on peut choisir une suite de poids $S=(S_k)_{k\geq 0} \subset (0,+\infty)$ pour les normes ${\cal C}^h$ de telle sorte à obtenir une estimation du type $\|P_q u\|_{r,h}\leq C\cdot \|u\|_{r,h}$, avec $C>0$ indépendante de la régularité $h\geq 0$ de la $q+1$-forme $u$. A condition de réstraindre opportunément le rayon $r>0$ on obtient un schéma de convergence rapide considérablement plus simple que  celui explique dans la preuve du théorème 1. En effet dans le cas on considération  on a pas besoin de réstraindre les rayons de la boule pendant les étapes du procédé itératif car on dispose de l'estimation précédente. Les détails de simplification et d'adaptation du procédé itératif  relatif à la preuve du théorème 1, au cas en examen sont laissé au lecteur. On obtient en conclusion une 
${\cal E}(\K) $-résolution locale, à partir d'une ${\cal E}(\K) $-résolution initiale, telle que les nouvelles matrices $\hat{\varphi}_j$ (ici on utilise les mêmes notations du théorème 1) soient toutes constantes. En particulier le fait que $\hat{\varphi}_1=Cst$ implique que ${\cal G}$ est un faisceau de ${\cal E}(\K)$-modules localement libre. La suite du théorème 2 dérive alors de résultats classiques bien connus du lecteur.
\\  
\\                                              
{\bf\large{Remerciements}}
\\
Je tiens enfin à remercier Jean-Pierre Demailly pour ses remarques et suggestions qui ont toujours constitué une aide précieuse pour la réalisation de ce travail.
%\vspace{\fill} 

\vspace{1cm}
\noindent
Nefton Pali
\\
Institut Fourier, UMR 5582, Université Joseph Fourier
\\
BP 74, 38402 St-Martin-d'Hères cedex, France
\\
E-mail: \textit{nefton.pali@ujf-grenoble.fr}
\end{document}